\documentclass[]{article}
\usepackage{amsmath}
\usepackage{amsfonts}
\usepackage{geometry}
\usepackage{graphicx}
\usepackage{tikz}
\usepackage{todonotes}

\newcommand{\subfigimg}[3][,]{%
  \setbox1=\hbox{\includegraphics[#1]{#3}}
  \leavevmode\rlap{\usebox1}
  \rlap{\hspace*{0pt}\raisebox{\dimexpr\ht1-0\baselineskip}{\bf
  \footnotesize #2}}
  \phantom{\usebox1}
}

\begin{document}

\title{Bugs on a Circle}
\author{Josh Briley \& Bryan Quaife}
\date{}
\maketitle

\begin{abstract}
We describe and analyze a generalization of the classic ``Four Bugs on a
  Square'' cyclic pursuit problem. Instead of allowing the bugs to
  spiral towards one another, we constrain $N$ bugs to the perimeter of
  the unit circle. Depending on their configuration, each bug moves
  either clockwise or counterclockwise with a constant angular speed, or
  remains stationary. Unlike the original problem where bugs always
  coalesce, this generalization produces three possible steady states:
  all bugs coalescing to a single point, clusters of bugs located at two
  antipodal points, or bugs entering a stable infinite chase cycle where
  they never meet. We analyze the stability of these steady states and
  calculate the probability that randomly initialized bugs reach each
  state. For $N \leq 4$, we derive exact analytical expressions for
  these probabilities. For larger values, we employ Monte Carlo
  simulations to estimate the probability of coalescing, finding it
  approximately follows an inverse square root relationship with the
  number of bugs. This generalization reveals rich dynamical behaviors
  that are absent in the classic problem. Our analysis provides insight
  into how restricting the bugs to the circle's perimeter fundamentally
  alters the long-term behavior of pursuing agents compared to
  unrestricted pursuit problems.
\end{abstract}

\section{Introduction}
\begin{figure}[h]
  \centering
  \subfigimg[height=0.4\textwidth]{(a)}{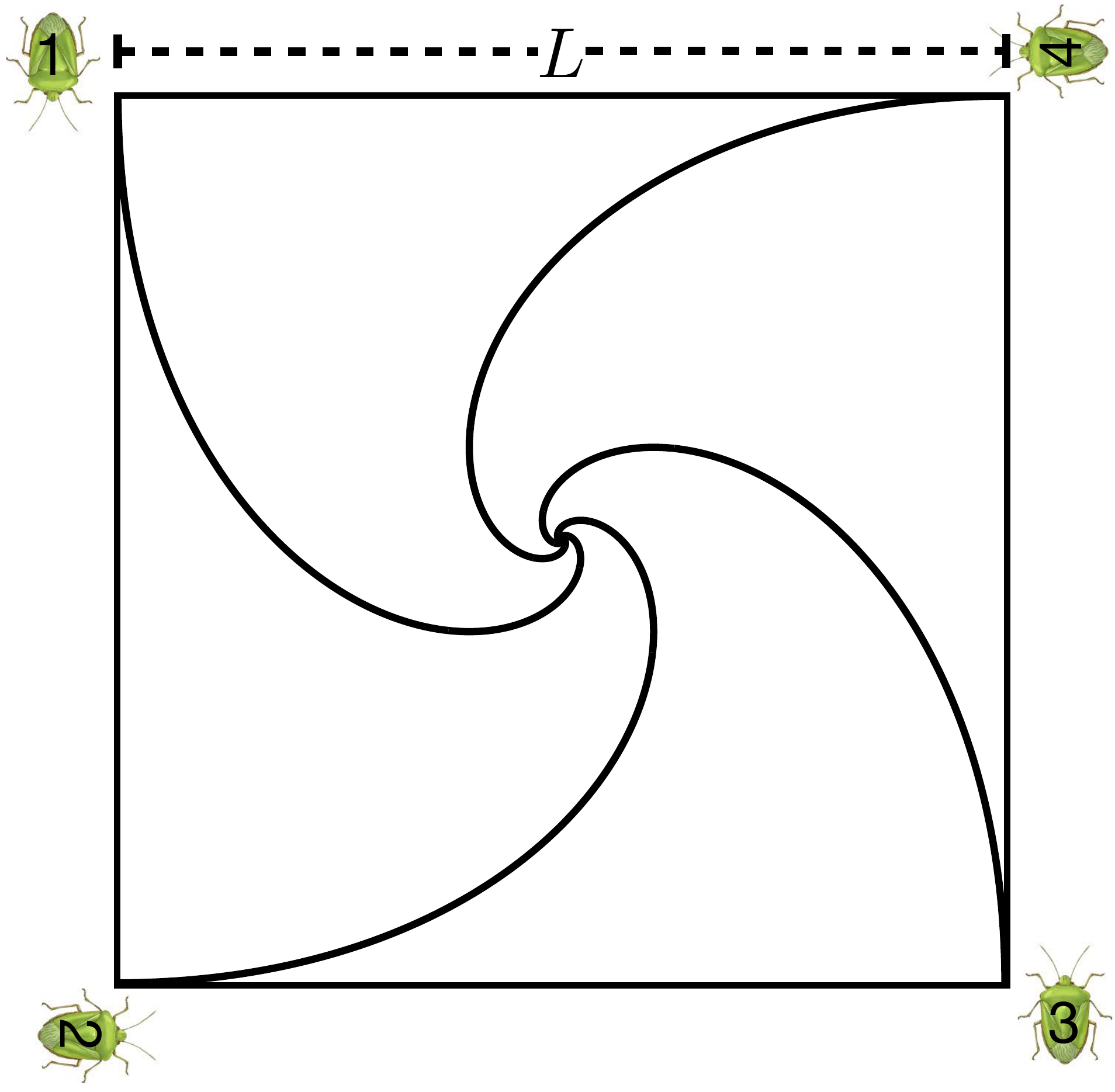} 
  \qquad \qquad
  \subfigimg[height=0.4\textwidth, trim=22cm 10cm 22cm 6cm,clip=true]{(b)}{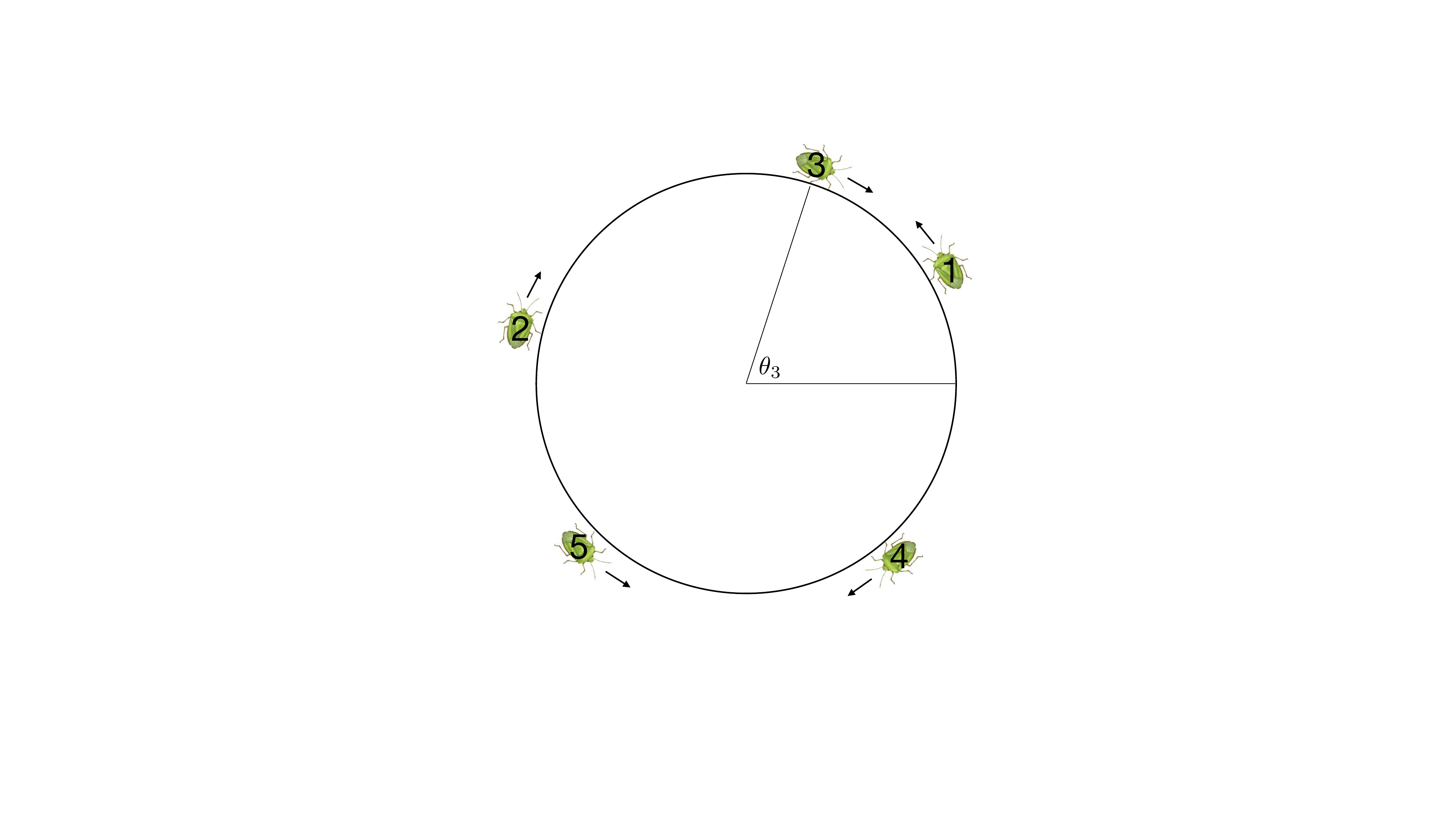}
  \caption{\label{fig:schematic} \em (a) Four bugs placed on the
  vertices of a square with sides of length $L$. The trajectories of the
  bugs are logarithmic spirals of length $L$ that meet at the middle of
  the square at time $L/V$. (b) Five bugs chasing each other around a
  circle. The angle $\theta_i$ is the polar angle of bug $i$. The arrows
  denote the direction the bugs travel in the current configuration.}
\end{figure}

In the classic four bugs on a square problem, four bugs are initially
placed at the corners of a square with length $L$, and each bug moves
towards its neighbor at a uniform speed $V$. This is done in a cyclic
fashion, so bug 1 chases bug 2, bug 2 chases bug 3, bug 3 chases bug 4,
and bug 4 chases bug 1. Early versions of this problem were described by
the popular science writer Martin Gardner, who called it one of {\em
Nine titillating puzzles}~\cite{gar1957}. The article asks questions
including: ``How far does each bug travel before they meet?'' In the
following month, Gardner showed that the bugs' trajectories form four
congruent logarithmic spirals of length $L$, that meet at the center of
the square at time $L/V$. Interestingly, the bugs rotate around the
center of the square infinitely many times. Gardner uses a geometric
argument to calculate the distance travelled, while Watton and Kydon use
Calculus~\cite{wat-kyd1969}. The four bugs problem has been discussed in
more recent news media---it appeared in the New York Times' Numberplay
on September 8, 2014. Figure~\ref{fig:schematic}(a) shows the
logarithmic spirals that each bug follows. The four bugs on a square
problem is an example of a general class of problems known as {\em
cyclic pursuit} which dates back to the nineteenth century. A detailed
historical description of the problem is provided by
Nahin~\cite{nah2012}.

Researchers have extended the model so that the bugs reach one of two
steady states: they either meet at a common point or they enter into a
circular cycle. Baryshnikov and Chen~\cite{bar-che2016} use the
terminology {\em rendez-vous} and {\em circular choreography} to mean
the case when bugs meet at a mutual point or when they enter into a
circular cycle, respectively. In this paper, we use the terminology {\em
coalesce} and {\em cycle}. Therefore, in the classic four bugs on a
square problem, the bugs coalesce. One obvious generalization is to
initialize the $N$ bugs on the vertices of a regular $N$-sided
polygon~\cite{gar1965}, and this results in the bugs coalescing.
However, by letting $N \rightarrow \infty$, the bugs remain on a circle,
resulting in them cycling. Other generalizations include the
introduction of non-linear or deviated steering laws that determine the
angular velocity of each bug. This can result in a cycle, and examples
of stable paths that the bugs will follow can be
circular~\cite{mukherjee2015deviated}, polygonal~\cite{mar-bro-fra2004},
or a figure-eight~\cite{bar-che2016}. Another way to produce cycling
behavior is to introduce time-dependent speed~\cite{bru-coh-efr1991,
kla-new1971}. Finally, stability of self-similar shapes that the bugs
follow when they coalesce has also been analyzed~\cite{ric2001}.

Another generalization of the classic problem involves initializing the bugs at positions other than the corners of a regular polygon. For example, four bugs can be
initialized at the four corners of a $2 \times 1$
rectangle~\cite{cha-lot-tre2011}. In this case, the bugs' locations are
always the vertices of a parallelogram, and that they meet in the middle
at time $1.5L/V$. The bugs coalesce at the middle of the initial
rectangle, and do so very rapidly---after each bug has completed one
rotation, they again form a 2-1 rectangle, but it has shrunk by a factor
of $10^{427907250}$. A broader generalization of cyclic pursuit are
evasion problems. These problems often have two groups of agents that
are trying to evade another group of agents. Seven Classic Evasion
Problems are described by Nahin~\cite{nah2012}. For example, in the
``Princess-and-Monster Game'', two agents are moving at discrete points
along the perimeter of a circle, but they are unaware of one another's
moves. This problem relates to ours, which we describe next, since the
agents are constrained to the perimeter of the circle.

In this work, we initialize $N$ bugs on the perimeter of the unit
circle, and require the bugs to only move in the tangential direction
with uniform speed. Interpreting each bug's location as a point on the
unit circle in the complex plane, the $N$ bugs at time $t$ are located
at $z_j(t) = e^{i\theta_j(t)}$, $\theta_j \in [0,2\pi)$, $j=1,\ldots,N$.
As in the case of a cyclic pursuit problem, each bug moves toward its
neighbor in a cyclic fashion, but the motion is restricted to the
angular direction at a constant angular velocity. We assume that each
bug can always see its neighboring bug across the circle, thereby
allowing it to know whether to move clockwise or counterclockwise. If
bug $i$ and bug $i+1$ meet, then they form a cluster and collectively
chase bug $i+2$. However, if bugs $i$ and $j$ meet, with $|i-j| \neq 1$,
then they simply pass through one another and continue to chase bugs
$i+1$ and $j+1$, respectively. Figure~\ref{fig:schematic}(b) shows a
configuration of $N=5$ bugs, and the direction that each bug is
currently traveling. Given this current configuration, bugs 1 and 2 will
eventually cluster and pursue bug 3, and bugs 4 and 5 will eventually
cluster and chase bug 1.

\begin{figure}[htp]
  \centering
  \subfigimg[width=0.31\textwidth,trim=14cm 6cm 14cm
  2cm,clip=true]{(a)}{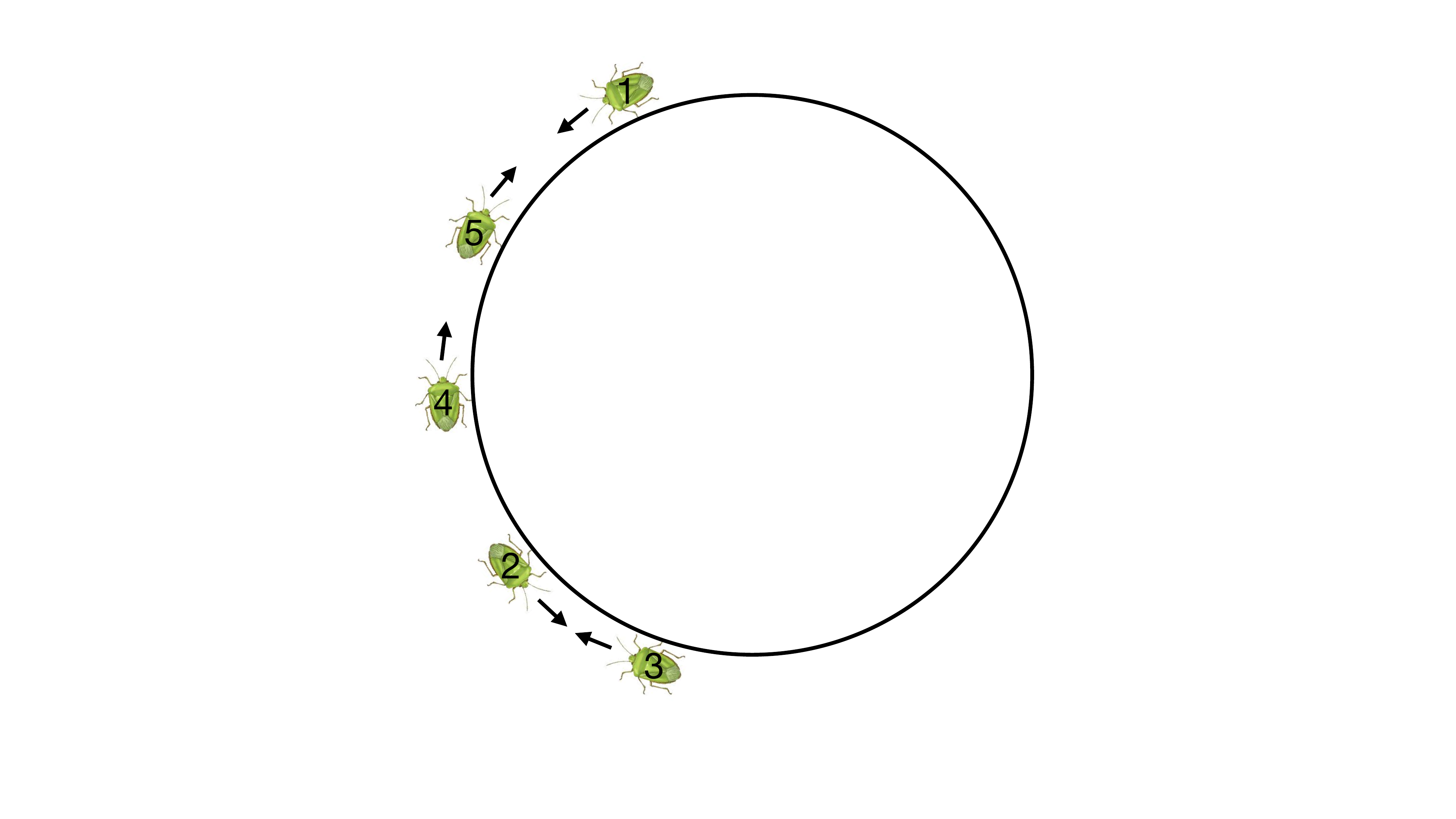}
  \quad
  \subfigimg[width=0.31\textwidth,trim=14cm 6cm 14cm
  2cm,clip=true]{(b)}{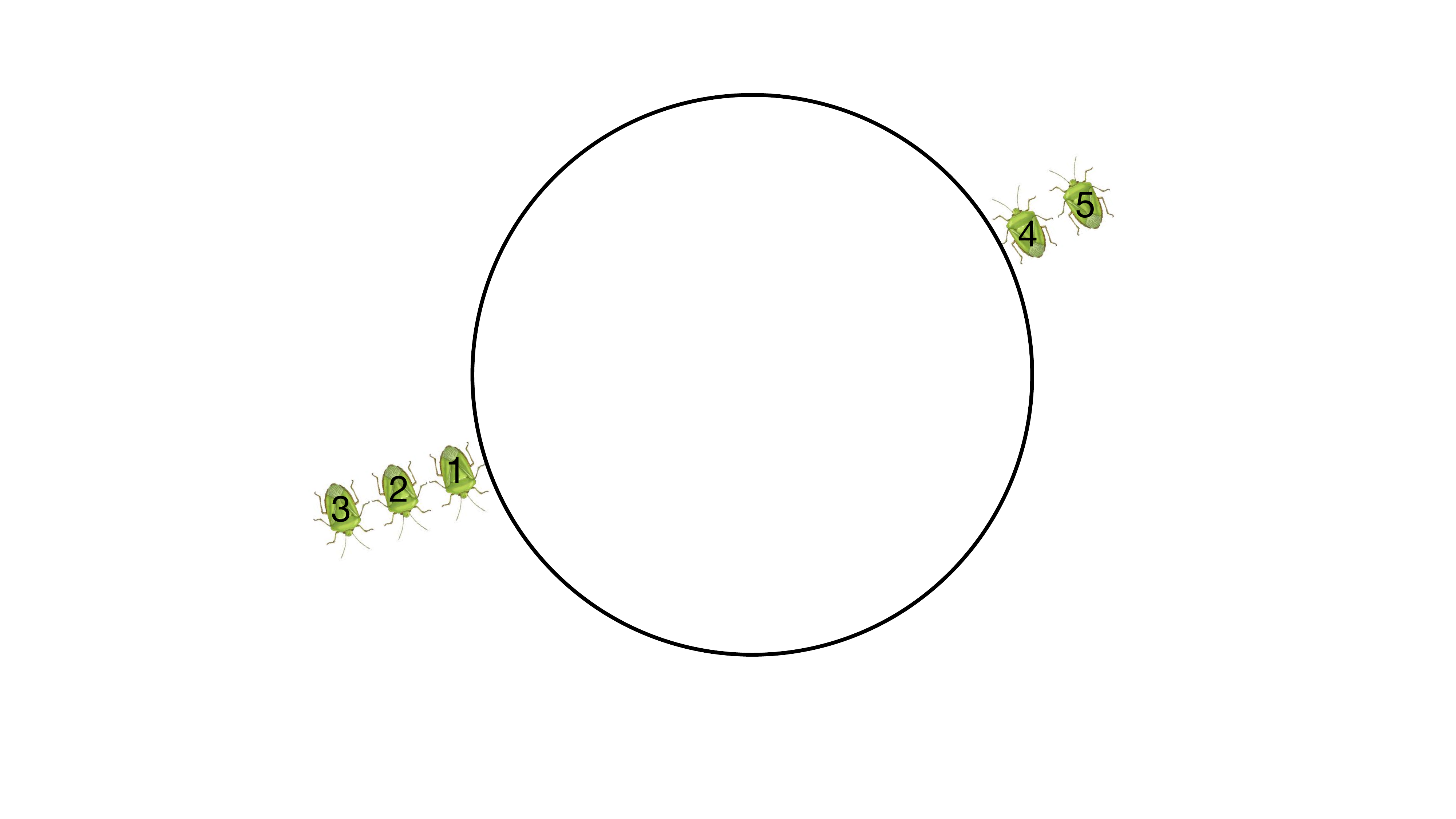}
  \quad
  \subfigimg[width=0.31\textwidth,trim=14cm 6cm 14cm
  2cm,clip=true]{(c)}{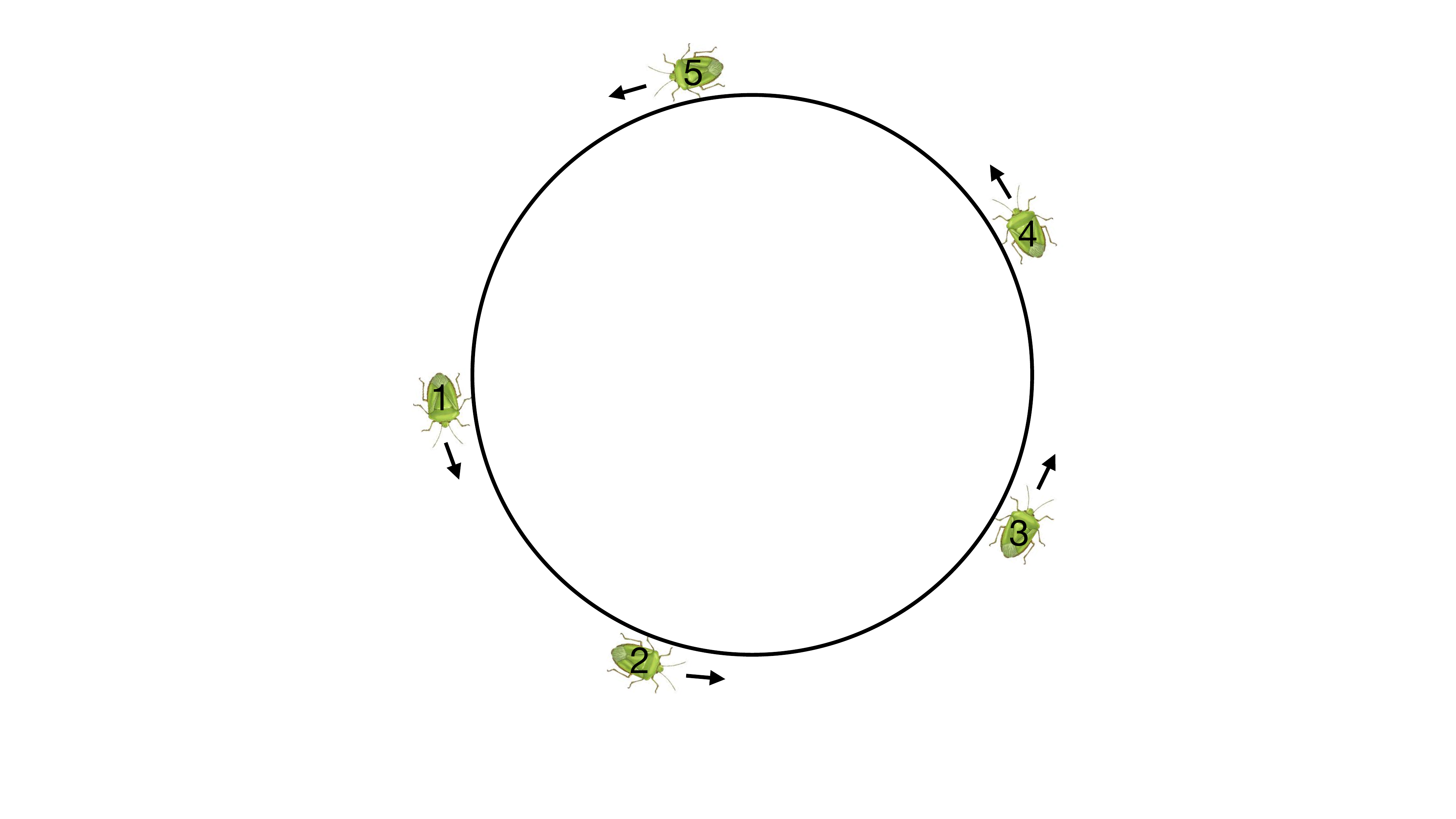}
  \caption{\label{fig:SSconfigs} \em Three possible steady states. (a)
  All the bugs are located on the same side of a diameter of the circle.
  The bugs will eventually coalesce at a single point which is a stable
  steady state. (b) Groups of bugs are located directly across from one
  another. This steady state is unstable. (c) Each bug is moving in a
  counterclockwise direction, resulting in a cycle. This steady state is
  stable.}
\end{figure}

Just like in the classic problem, the bugs can coalesce to a single
point. The inevitable coalescence of all bugs to a single point can be
detected if, at any point in time, all bugs are located on one side of a
diameter of the circle (Figure~\ref{fig:SSconfigs}(a)). However, a major
difference introduced by this generalization is the existence of
additional steady states without requiring any other generalizations of
the problem. First, two groups of bugs can be located at antipodal
points on the circle (Figure~\ref{fig:SSconfigs}(b)). This configuration
results in the bugs not moving, but it is an unstable steady state. Any
perturbation to the bugs results in the system tending to another steady
state. Second, the bugs can be distributed around the circle so that
they are all moving in the same direction. This results in a cycle where
the bugs continuously move, but never reach one another
(Figure~\ref{fig:SSconfigs}(c)). This configuration can be easily
detected if at any point in time, all the bugs are moving in the same
direction.

In this paper, we analyze these steady states by considering their
stability, and calculating probabilities that randomly initialized bugs
reach one of these states. Section~\ref{sec:governing} introduces the
notation and governing equations. Section~\ref{sec:small} analyzes
stability and probabilities analytically for the case $N \leq 4$. Then
Section~\ref{sec:large} describes a numerical method to simulate the
bugs' dynamics and uses Monte Carlo methods to study the stability and
probabilities for large numbers of bugs.

\section{Governing Equations and Notation}
\label{sec:governing}
Given bugs located at locations $z_j(t) = \exp(i\theta_j(t))$,
$j=1,\ldots,N$, bug $j$ is governed by the equation
\begin{align}
  \label{eqn:dynamics}
  \frac{d}{dt}\theta_j(t) = 
  \begin{cases}
    +1, & \mod(\theta_{j+1} - \theta_{j},2\pi) \in (0,\pi), \\
    -1, & \mod(\theta_{j+1} - \theta_{j},2\pi) \in (\pi,2\pi), \\
    0,  & \mod(\theta_{j+1} - \theta_{j},2\pi) \in \{0,\pi\}.
  \end{cases}
\end{align}
Since we assume cyclic pursuit, we associate bug $(N+1)$ with bug $1$.
It is convenient to maintain each bug's angle within the interval
$[0,2\pi)$. Therefore, if a bug passes through the point $(1,0)$
(i.e.~$\theta_j=0$), then we either subtract or add $2\pi$ to the angle.
For simplicity, we omit the $\mod$ notation throughout, but note that we
consider all angles to be defined in the principal branch $[0,2\pi)$.

Our model is closely related to the Kuramoto
model~\cite{kuramoto1984chemical}
\begin{align}
  \label{eqn:kuramoto}
  \dot{\theta}_j = \omega_j + \sum_{k \neq j} K_{jk} 
    \Phi(\theta_j - \theta_k),
\end{align}
which is used to describe synchronization of chemical and biological
oscillators. Equation~\eqref{eqn:dynamics} reduces to
equation~\eqref{eqn:kuramoto} when $\omega_j = 0$, $\Phi(z) = z$, and 
\begin{align}
  K_{jk} = \begin{cases}
    1, & j = k+1, \\
    0, & \text{otherwise.}
  \end{cases}
\end{align}
While Kuramoto models can exhibit chaos for sufficiently large $N$,
chaos is not expected for solutions to~\eqref{eqn:dynamics}. Although,
as we will see, it can take a long time before the system converges to a
steady state.

To study the stability of equation~\eqref{eqn:dynamics}, it is sometimes
convenient to decrease the number of unknowns from $N$ to $N-1$ by
arbitrarily rotating the whole system at each time so that bug 1 is
always located at $\theta_1 = 0$. We introduce the variables
\begin{align}
  \omega_{j} = \theta_{j+1} - \theta_j, \quad j=1,\ldots,N-1,
\end{align}
which measures the angle between neighboring bugs. The governing
equations for the angle $\omega_j$ depend on the relative location of
$\theta_{j+1}$ and $\theta_{j}$, and the relative location of
$\theta_{j}$ and $\theta_{j-1}$. In particular
\begin{align}
  \label{eqn:omegaDynamics}
  \frac{d}{dt}\omega_j(t) = 
  \begin{cases}
    +2, & \omega_{j+1} \in (0,\pi), \: \text{and} \:
          \omega_{j} \in (\pi,2\pi), \\
    -2, & \omega_{j+1} \in (\pi,2\pi), \: \text{and} \:
          \omega_{j} \in (0,\pi), \\
    0, & \text{otherwise}.
  \end{cases}
\end{align}
The case $\dot{\omega}_j(t)=0$ occurs when both bugs $j$ and $j+1$ are
moving in the same direction.

In addition to studying the problem's dynamics, we consider the problem
from a statistical point of view. Under a random initialization of the
$N$ bug locations, the system evolves to one of three states: $X_N =
\texttt{coalesce}$ (Figure~\ref{fig:SSconfigs}(a)), $X_N =
\texttt{groups}$ (Figure~\ref{fig:SSconfigs}(b)), or $X_N =
\texttt{cycle}$ (Figure~\ref{fig:SSconfigs}(c)). As mentioned earlier,
the \texttt{groups} steady state is unstable, meaning that $N$ randomly
initialized bugs reach this state is with probability 0. Therefore,
$X_N$ is a binary random variable that takes on one of the two values
\begin{align}
  X_N = \{\texttt{coalesce},\texttt{cycle}\}.
\end{align}
We define $p_N$ as the probability that $N$ randomly initialized bugs
reach a \texttt{cycle} state, and then $1-p_N$ is the probability that
$N$ randomly initialized bugs reach a \texttt{coalesce} state. One of
the main goals of this paper is to calculate $p_N$ for small $N$ and
estimate $p_N$ using Monte Carlo methods for large $N$.

\section{Analytic Calculations for Small $N$}
\label{sec:small}
For sufficiently small $N$, the stability of the steady states, and
probabilities that randomly initialized bugs reach such a state can be
calculated analytically. In this section, we consider the cases $N = 2$,
$3$, and $4$.

\subsection{Two Bug Case}
Although trivial, it is insightful to first consider the two-bug case.
The governing equation for $\omega = \theta_2 - \theta_1$ is
\begin{align}
  \frac{d\omega}{dt} = f(\omega) = 
  \begin{cases}
    -2, & \omega \in (0,\pi), \\
    +2, & \omega \in (\pi,2\pi), \\ 
     0, & \text{otherwise.}
  \end{cases}
\end{align}
The two fixed points are clearly $\omega = 0$ and $\omega = \pi$ which
correspond to the cases $\theta_1 = \theta_2$ and $\theta_2 - \theta_1 =
\pi$, respectively. We analyze the stability of each of these fixed
points by considering small perturbations. Letting $\epsilon > 0$, we
have that $f(0 + \epsilon) < 0$ and $f(0 - \epsilon)
> 0$, showing that the system will return towards $0$ under a
small perturbation, rendering $\omega=0$ a stable steady state.
Conversely, $f(\pi + \epsilon) > 0$ and $f(\pi - \epsilon) < 0$, showing
that the system will diverge from $\pi$ under a small perturbation,
rendering $\omega = \pi$ an unstable steady state. Alternatively, we can
analyze the stability of each of these steady states by noting that 
\begin{align}
  f'(\omega) = 4\delta(\omega - \pi) - 4\delta(\omega - 0),
\end{align}
where we have taken the derivative in the weak sense, and $\delta$ is
the Dirac delta function. With a slight abuse of notation, we see that
$f'(\pi) > 0$ and $f'(0) < 0$, showing that the $\omega = \pi$ is an
unstable fixed point while $\omega = 0$ is a stable fixed point. Even
without this analysis, the stability is intuitively clear since the
case $\omega = \pi$ corresponds to bugs being antipodal to one another,
and any small perturbation would result in the bugs moving towards the
stable fixed point $\omega = 0$, the case where the bugs are at the same
location. In conclusion, given an initial condition of the bug
locations, we have that the bugs will always coalesce, meaning that
\begin{align}
  p_2 = 0.
\end{align}

\subsection{Three Bug Case}
The $N=3$ case introduces interesting behaviors that are absent in the
$N=2$ case. We first write the governing equations in terms of $\omega_1
= \theta_2 - \theta_1$ and $\omega_2 = \theta_3 - \theta_2$. To see that
there is a stable cycle steady state, meaning that $p_3 > 0$, suppose
that the bugs are initialized at $\theta_1 = 0$, $\theta_2 = 2\pi/3$ and
$\theta_3 = 4\pi/3$. In this case, $\omega_1 = \omega_2 = 2\pi/3$. Since
$\omega_1 \in (0,\pi)$ and $\omega_2 \in (0,\pi)$,
equation~\eqref{eqn:omegaDynamics} shows that $\dot{\omega}_1 =
\dot{\omega}_2 = 0$, meaning that the pairwise distances do not change.
In this configuration, all three bugs rotate counterclockwise at a
uniform speed. Furthermore, under a sufficiently small perturbation of
the bug locations, we will still have $\dot{\omega}_1 = \dot{\omega}_2 =
0$, meaning that this steady state is stable. Finally, to confirm that
the bugs can still \texttt{coalesce}, meaning that $p_3 < 1$, recall
that if all bugs are located on one side of any diameter of the circle,
then they will \texttt{coalesce}.

The $N = 3$ case is sufficiently simple that $p_3$ can be calculated
analytically. We calculate the probability that the bugs enter the
\texttt{cycle} steady state and rotate counterclockwise. We can perform
this calculation in terms of $\omega_1$ and $\omega_2$, but it is
convenient to use the variables $\theta_1$, $\theta_2$, and $\theta_3$.
The system can reach this steady state only if the bugs are initialized
so that they are all moving counterclockwise at $t=0$. Furthermore,
because the system is invariant to rotation, we can assume that
$\theta_1 = 0$. Therefore, 
\begin{align}
  P(X_3 = \texttt{cycle} \text{ and } \texttt{counterclockwise}) &= 
  P((\theta_2 \in (0,\pi)) \cap (\theta_3 \in (\pi,\theta_2 + \pi))) \\
  &= \left(\frac{1}{2\pi}\right)^2 \int_0^\pi \int_\pi^{\theta_2 + \pi} 
    d\theta_3 d\theta_2 = \frac{1}{8}.
\end{align} 
Doubling the result to account for clockwise rotation, we have that
\begin{align}
  \label{eqn:3BugsProbs}
  p_3 = \frac{1}{4}.
\end{align}

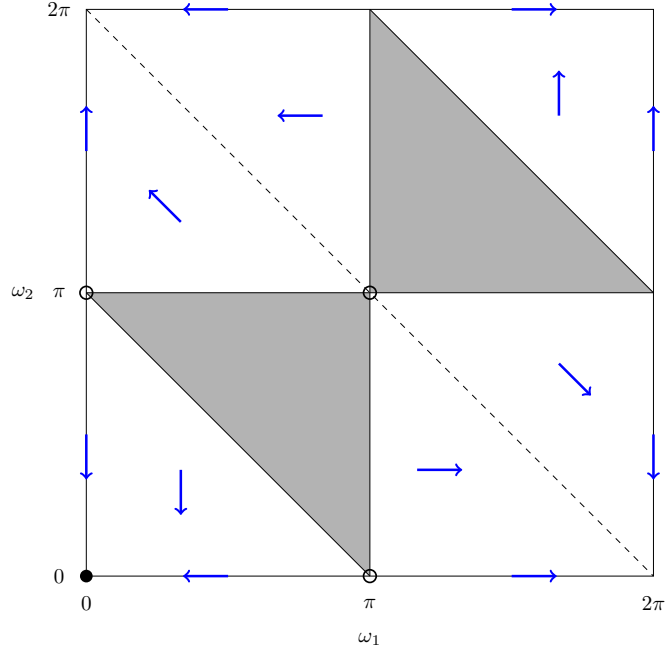
\begin{figure}[htp]
  \centering
  \scalebox{0.8}{\begin{tikzpicture}[scale=1.5]

  \fill[black!30] (0, pi) -- (pi, pi) -- (pi, 0) -- cycle;
  \fill[black!30] (pi, pi) -- (pi, 2*pi) -- (2*pi, pi) -- cycle;

  \draw[thin] (0,0) rectangle ({2*pi},{2*pi});

  \node at (pi,-0.7) {$\omega_1$};
  \node at (-0.7,pi) {$\omega_2$};


  \node at (0,-0.3) {$0$};
  \node at (pi,-0.3) {$\pi$};
  \node at (2*pi,-0.3) {$2\pi$};
  \node at (-0.3,0) {$0$};
  \node at (-0.3,pi) {$\pi$};
  \node at (-0.3,2*pi) {$2\pi$};

  \draw[thin] (0*pi, 1*pi) -- (1*pi, 0*pi);
  \draw[thin,dashed] (0*pi, 2*pi) -- (2*pi, 0*pi);
  \draw[thin] (1*pi, 2*pi) -- (2*pi, 1*pi);
  \draw[thin] (0*pi, 1*pi) -- (2*pi, 1*pi);
  \draw[thin] (1*pi, 0*pi) -- (1*pi, 2*pi);

  \draw[->, blue, very thick] (pi/3, 5*pi/4) -- ++(-0.7071/2, +0.7071/2);  
  \draw[->, blue, very thick] (5*pi/6, 13*pi/8) -- ++(-1/2, +0);  
  \draw[->, blue, very thick] (5*pi/3, 13*pi/8) -- ++(+0, +1/2);  

  \draw[->, blue, very thick] (pi/3, 3*pi/8) -- ++(+0, -1/2);  
  \draw[->, blue, very thick] (7*pi/6, 3*pi/8) -- ++(+1/2, +0);  
  \draw[->, blue, very thick] (5*pi/3, 3*pi/4) -- ++(+0.7071/2, -0.7071/2);  
  
  \draw[->, blue, very thick] (pi/2,0) -- ++(-1/2, +0);  
  \draw[->, blue, very thick] (3*pi/2,0) -- ++(+1/2, +0);  
  \draw[->, blue, very thick] (2*pi,pi/2) -- ++(+0, -1/2);  
  \draw[->, blue, very thick] (2*pi,3*pi/2) -- ++(+0, +1/2);  
  \draw[->, blue, very thick] (3*pi/2,2*pi) -- ++(+1/2, +0);  
  \draw[->, blue, very thick] (pi/2,2*pi) -- ++(-1/2, +0);  
  \draw[->, blue, very thick] (0,pi/2) -- ++(+0, -1/2);  
  \draw[->, blue, very thick] (0,3*pi/2) -- ++(+0, +1/2);  

  \fill (0,0) circle[radius=2pt];

  \draw[thick] (pi,0) circle[radius=2pt];
  \draw[thick] (0,pi) circle[radius=2pt];
  \draw[thick] (pi,pi) circle[radius=2pt];

\end{tikzpicture}}
  \caption{\label{fig:phaseSpace} \em The phase portrait of the $N = 3$
  case. The blue arrows show the trajectory of $(\omega_1,\omega_2)$.
  The bugs \texttt{cycle} if $(\omega_1,\omega_2)$ is located in the
  gray regions. The stable \texttt{coalesce} point (filled black point)
  is $(0,0)$. The unstable \texttt{groups} points (open black points)
  are $(0,\pi)$, $(\pi,\pi)$, and $(\pi,0)$. The dashed black line
  separates regions where both $\omega_1$ and $\omega_2$ are changing,
  versus where $\omega_1$ is changing while $\omega_2$ is not. Note that
  the figure is doubly-periodic.}
\end{figure}

We gain further insight by using a geometric argument by visualizing a
phase portrait. Figure~\ref{fig:phaseSpace} shows all possible
configurations of the system in terms of $(\omega_1,\omega_2) \in
[0,2\pi)^2$. The vectors show the trajectory that $\omega_1$ and
$\omega_2$ will follow towards the steady states. If
$(\omega_1,\omega_2)$ are located in one of the gray triangles, the bugs
\texttt{cycle} either clockwise (upper right triangle) or
counterclockwise (lower left triangle). The unstable \texttt{groups}
steady states are located at the open black circles, and the stable
\texttt{coalesce} steady state is at the closed black circle. Points
between the open circles ($\omega_1 = \pi$ and $\omega_2 = \pi$)
correspond to an unstable \texttt{cycle} steady state.  Finally, the
bugs reach a stable \texttt{coalesce} steady state if
$(\omega_1,\omega_2)$ is located interior to a white triangle, if
$\omega_1 = 0$, or if $\omega_2 = 0$. We can easily verify that the
probability in equation~\eqref{eqn:3BugsProbs} is correct by noting that
the area of the gray region comprises one quarter of the total area of
the square.

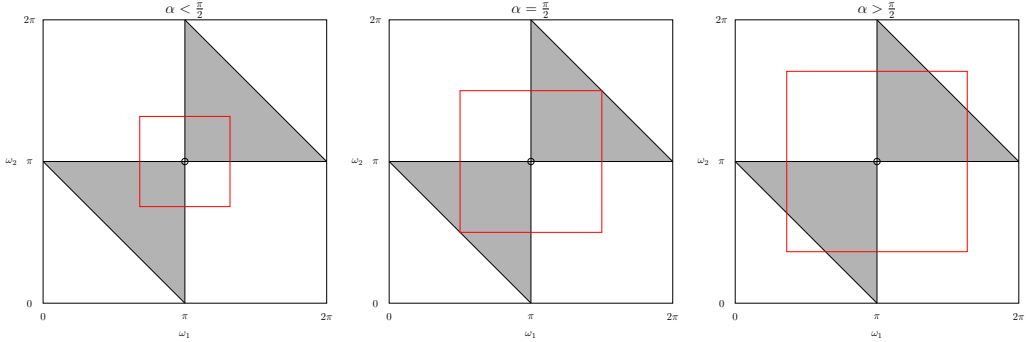
\begin{figure}[htp]
  \centering
  \scalebox{0.4}{\begin{tikzpicture}[scale=1.5]

  \fill[black!30] (0, pi) -- (pi, pi) -- (pi, 0) -- cycle;
  \fill[black!30] (pi, pi) -- (pi, 2*pi) -- (2*pi, pi) -- cycle;

  \draw[thin] (0,0) rectangle ({2*pi},{2*pi});

  \node at (pi,-0.7) {$\omega_1$};
  \node at (-0.7,pi) {$\omega_2$};

  \node at (0,-0.3) {$0$};
  \node at (pi,-0.3) {$\pi$};
  \node at (2*pi,-0.3) {$2\pi$};
  \node at (-0.3,0) {$0$};
  \node at (-0.3,pi) {$\pi$};
  \node at (-0.3,2*pi) {$2\pi$};

  \draw[thin] (0*pi, 1*pi) -- (1*pi, 0*pi);
  \draw[thin] (1*pi, 2*pi) -- (2*pi, 1*pi);
  \draw[thin] (0*pi, 1*pi) -- (2*pi, 1*pi);
  \draw[thin] (1*pi, 0*pi) -- (1*pi, 2*pi);

  \draw[thick] (pi,pi) circle[radius=2pt];

  \draw[thin,red] (pi-1, pi-1) -- (pi-1, pi+1) -- 
      (pi+1,pi+1) -- (pi+1,pi-1) -- cycle;
  \node at (pi,6.5) {\Large $\alpha < \frac{\pi}{2}$};

\end{tikzpicture}}
  \scalebox{0.4}{\begin{tikzpicture}[scale=1.5]

  \fill[black!30] (0, pi) -- (pi, pi) -- (pi, 0) -- cycle;
  \fill[black!30] (pi, pi) -- (pi, 2*pi) -- (2*pi, pi) -- cycle;

  \draw[thin] (0,0) rectangle ({2*pi},{2*pi});

  \node at (pi,-0.7) {$\omega_1$};
  \node at (-0.7,pi) {$\omega_2$};

  \node at (0,-0.3) {$0$};
  \node at (pi,-0.3) {$\pi$};
  \node at (2*pi,-0.3) {$2\pi$};
  \node at (-0.3,0) {$0$};
  \node at (-0.3,pi) {$\pi$};
  \node at (-0.3,2*pi) {$2\pi$};

  \draw[thin] (0*pi, 1*pi) -- (1*pi, 0*pi);
  \draw[thin] (1*pi, 2*pi) -- (2*pi, 1*pi);
  \draw[thin] (0*pi, 1*pi) -- (2*pi, 1*pi);
  \draw[thin] (1*pi, 0*pi) -- (1*pi, 2*pi);

  \draw[thick] (pi,pi) circle[radius=2pt];

  \draw[thin,red] (pi-pi/2, pi-pi/2) -- (pi-pi/2, pi+pi/2) -- 
      (pi+pi/2,pi+pi/2) -- (pi+pi/2,pi-pi/2) -- cycle;
  \node at (pi,6.5) {\Large $\alpha = \frac{\pi}{2}$};

\end{tikzpicture}}
  \scalebox{0.4}{\begin{tikzpicture}[scale=1.5]

  \fill[black!30] (0, pi) -- (pi, pi) -- (pi, 0) -- cycle;
  \fill[black!30] (pi, pi) -- (pi, 2*pi) -- (2*pi, pi) -- cycle;

  \draw[thin] (0,0) rectangle ({2*pi},{2*pi});

  \node at (pi,-0.7) {$\omega_1$};
  \node at (-0.7,pi) {$\omega_2$};

  \node at (0,-0.3) {$0$};
  \node at (pi,-0.3) {$\pi$};
  \node at (2*pi,-0.3) {$2\pi$};
  \node at (-0.3,0) {$0$};
  \node at (-0.3,pi) {$\pi$};
  \node at (-0.3,2*pi) {$2\pi$};

  \draw[thin] (0*pi, 1*pi) -- (1*pi, 0*pi);
  \draw[thin] (1*pi, 2*pi) -- (2*pi, 1*pi);
  \draw[thin] (0*pi, 1*pi) -- (2*pi, 1*pi);
  \draw[thin] (1*pi, 0*pi) -- (1*pi, 2*pi);

  \draw[thick] (pi,pi) circle[radius=2pt];

  \draw[thin,red] (pi-2, pi-2) -- (pi-2, pi+2) -- 
      (pi+2,pi+2) -- (pi+2,pi-2) -- cycle;
  \node at (pi,6.5) {\Large $\alpha > \frac{\pi}{2}$};

\end{tikzpicture}}
  \caption{\label{fig:Stability3Bug} \em While the fixed point of
  $(\omega_1,\omega_2) = (\pi,\pi)$ is unstable, small perturbations can
  result in either a \texttt{coalesce} or \texttt{cycle} steady state.
  The probabilities of entering the \texttt{cycle} steady state is the
  area of the gray region inside the red square relative to the area of
  the red square. The probability is given in
  equation~\eqref{eqn:stabilityProb}.}
\end{figure}

We can use the phase portrait to analyze the instability of the
\texttt{groups} steady state (open circles). All three of these points
have the same stability behavior, so we focus on $(\omega_1,\omega_2) =
(\pi,\pi)$, corresponding to bugs initialized at $\theta_1 = 0$,
$\theta_2 = \pi$, and $\theta_3 = 0$. Under a perturbation, we see that
the system may enter into the gray region, in which case the bugs will
\texttt{cycle}, or into the white region, in which case the bugs will
\texttt{coalesce}. Suppose we allow both $\omega_1$ and $\omega_2$ to be
perturbed by $0 < \Delta \omega < \alpha$
(Figure~\ref{fig:Stability3Bug}). If $\alpha \leq \pi/2$, then the
probability that the perturbed system goes to a \texttt{cycle} steady
state is 0.5. If $\alpha > \pi/2$, then the probability that the
perturbed system enters a \texttt{cycle} steady state is the area of the
gray region inside the red square. A geometric calculation shows that
\begin{align}
  \label{eqn:stabilityProb}
  p_3^{\texttt{stab}}(\alpha):=
  P(X_3 = \texttt{cycle} \: | \: \Delta \omega < \alpha) = 
  \begin{cases}
    \frac{1}{2}, & \text{if } \alpha \leq \frac{\pi}{2}, \\
    \frac{1}{4}\left(-\left(\frac{\pi}{\alpha}\right)^2 + 
    4\left(\frac{\pi}{\alpha}\right) - 2\right), 
    & \text{if } \alpha > \frac{\pi}{2}.
  \end{cases}
\end{align}
This result is used to validate our Monte Carlo approach in
Section~\ref{sec:large}.

\subsection{Four Bug Case}
\begin{figure}[h]
  \centering
  \subfigimg[width=0.48\textwidth, trim=15cm 7cm 15cm 3cm,clip=true]{(a)}{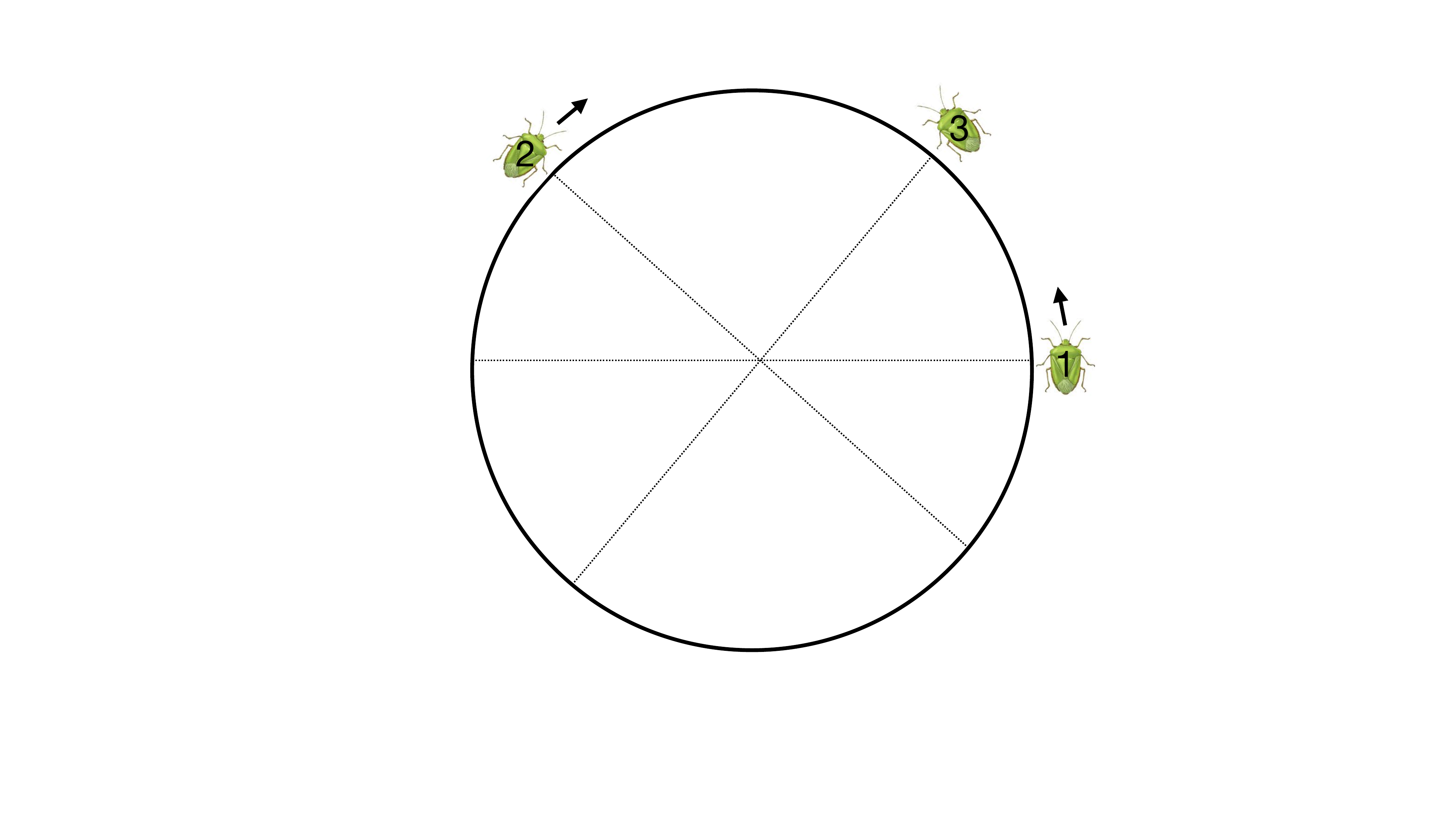}
  \subfigimg[width=0.48\textwidth, trim=15cm 7cm 15cm 3cm,clip=true]{(b)}{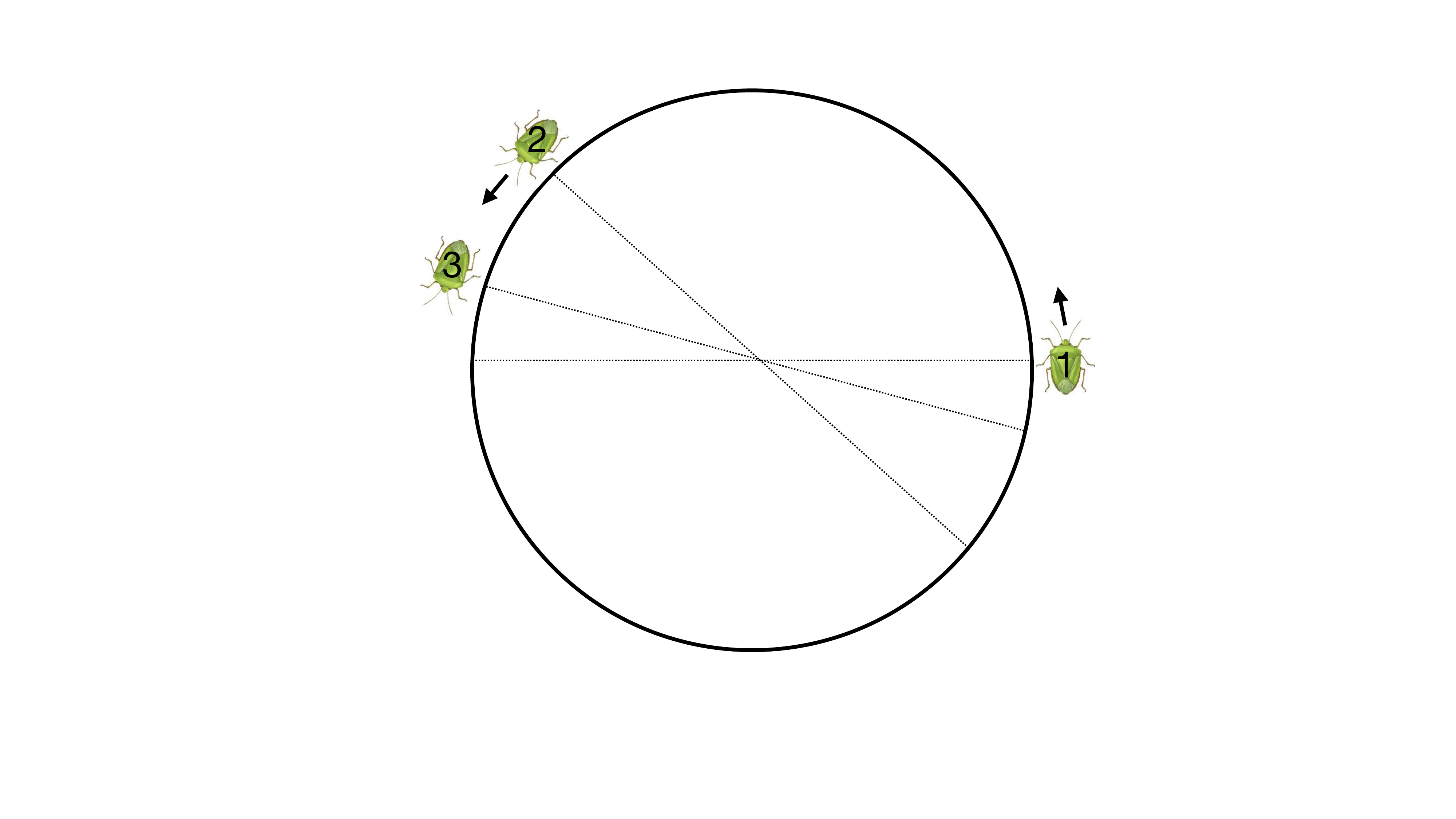}
  \subfigimg[width=0.48\textwidth, trim=15cm 7cm 15cm 3cm,clip=true]{(c)}{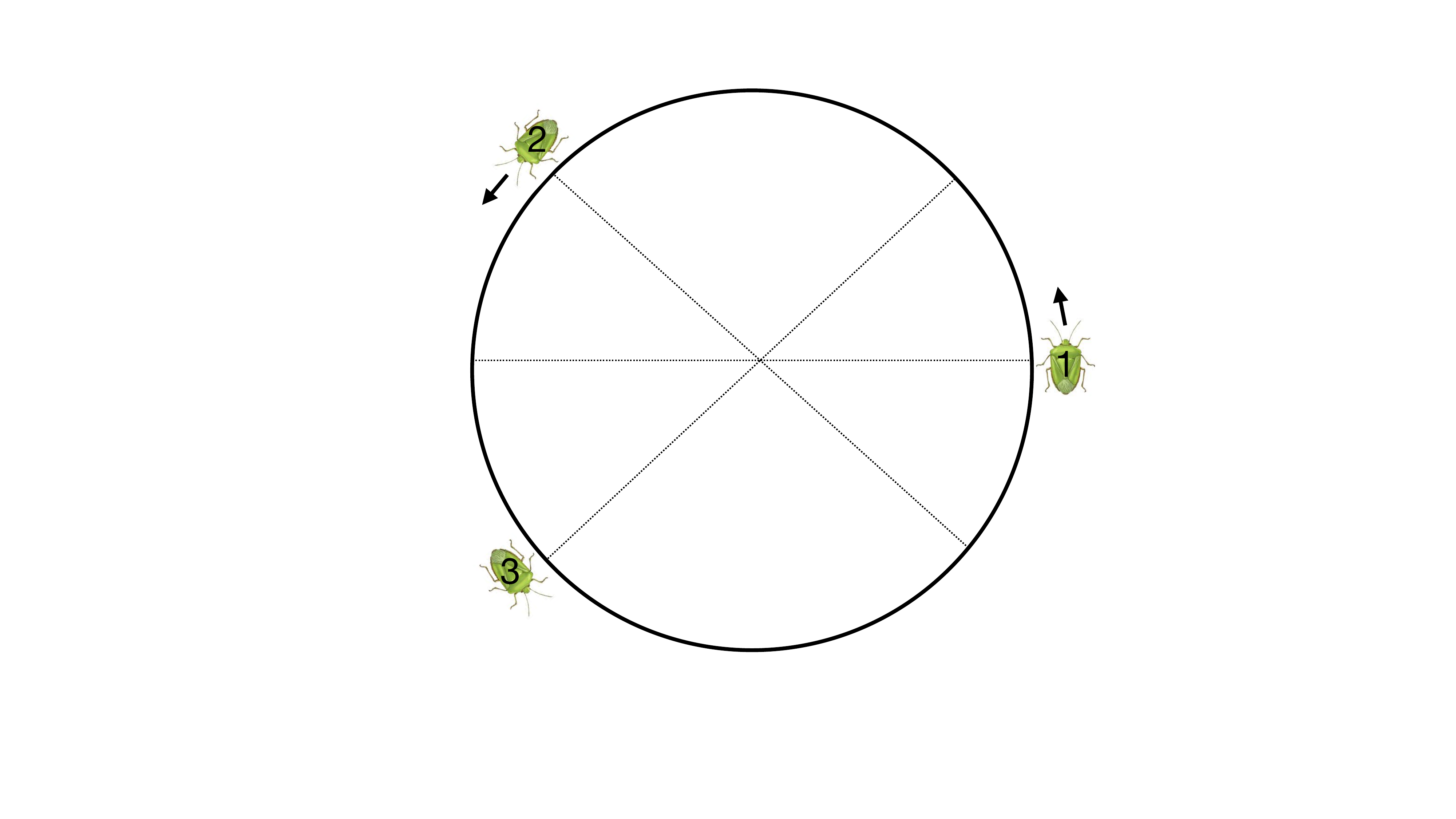}
  \subfigimg[width=0.48\textwidth, trim=15cm 7cm 15cm 3cm,clip=true]{(d)}{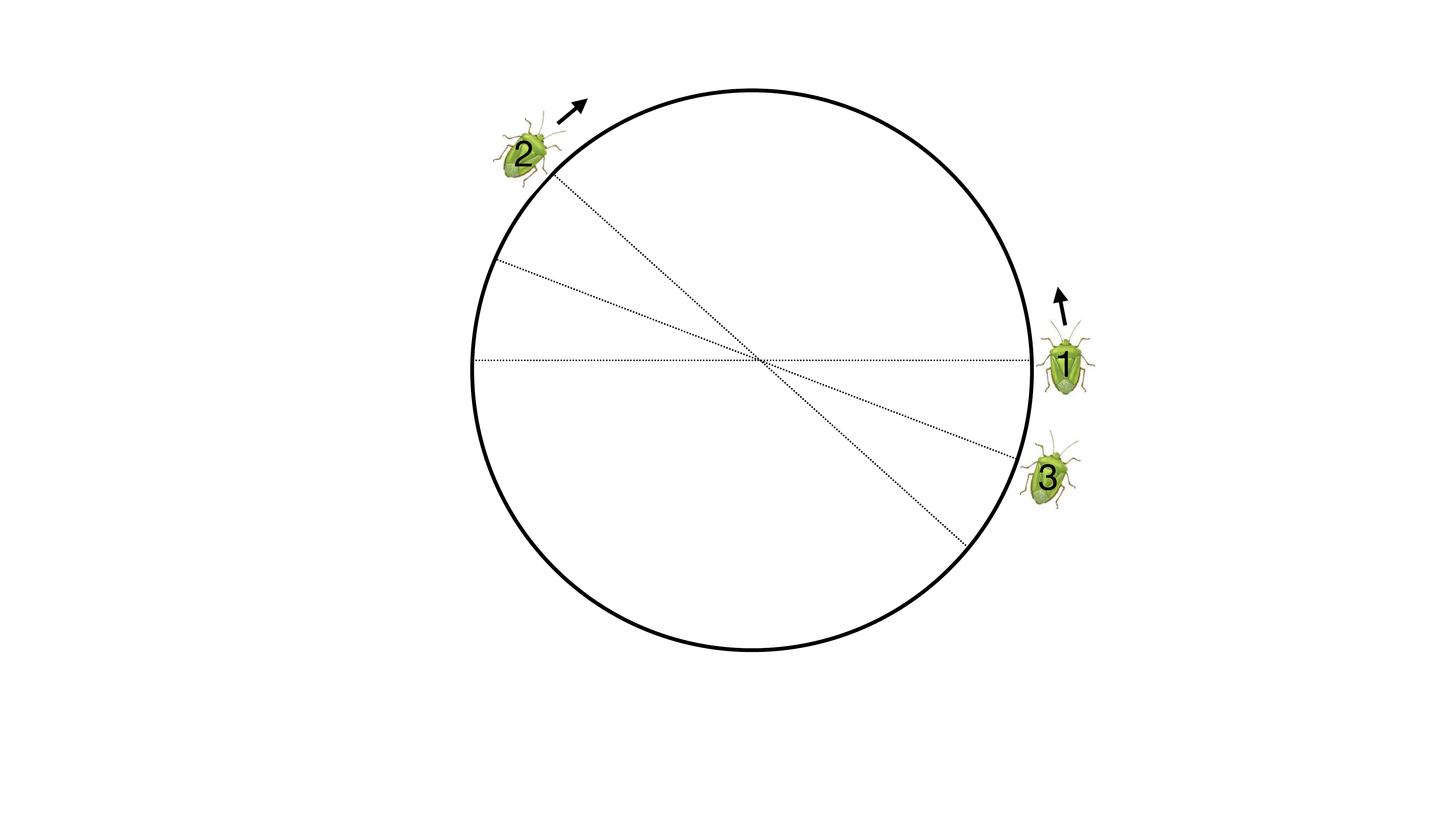}
  \caption{\label{fig:FourBugConfigs} \em Four possible initial
  conditions of bugs 1, 2, and 3 of a four bug system. We assume
  $\theta_1(0) = 0$ and $\theta_2(0) \in (0,\pi)$. Bug 3 is located in
  the interval (a) $(0,\theta_2)$, (b) $(\theta_2,\pi)$, (c) $
  (\pi,\theta_2 + \pi)$, (d) $(\theta_2 + \pi,2\pi)$.}
\end{figure}

We calculate $p_4$ as a final analytical case before resorting to Monte
Carlo methods. As before, we let $\theta_1(0) = 0$. We also make use of
symmetry by first assuming that $\theta_2 \in (0,\pi)$, and then
doubling the result to calculate $p_4$. To determine if a choice for
$\theta_3$ and $\theta_4$ will result in a \texttt{cycle} or a
\texttt{coalesce} steady state, we consider four regions of the
circle that are defined by the choice of $\theta_2$. The regions are
$(0,\theta_2)$, $(\theta_2,\pi)$, $(\pi,\theta_2 + \pi)$, and
$(\theta_2+\pi,2\pi)$. Figure~\ref{fig:FourBugConfigs} shows these four
possible regions that bug 3 can be initialized.

Each possible value for $\theta_3(0)$ defines six intervals. The
interval that contains $\theta_4(0)$ determines whether the bugs will
reach a \texttt{cycle} or \texttt{coalesce} steady state. It is obvious
in several of the cases if the bugs will \texttt{cycle} or
\texttt{coalesce}. For example, with the configuration in
Figure~\ref{fig:FourBugConfigs}(a), if bug 4 is located in
$(0,\theta_3)$, then all bugs are initially on the upper hemisphere of
the circle, so they will \texttt{coalesce}. In contrast, in
Figure~\ref{fig:FourBugConfigs}(c), if bug 4 is located in
$(\theta_3,\theta_2 + \pi)$, then the four bugs will \texttt{cycle}. In
each of the four cases for the initial placement of $\theta_3$, and each
of the six cases for the initial placement of $\theta_4$, the bugs will
eventually \texttt{cycle} if
\begin{itemize}
  \item $\theta_3 \in (0,\theta_2)$ and $\theta_4 \in (\pi,\theta_3 +
    \pi)$;

  \item $\theta_3 \in (\theta_2,\pi)$ and $\theta_4 \in (\pi,\theta_3 +
    \pi)$;

  \item $\theta_3 \in (\pi,\theta_2 + \pi)$ and $\theta_4 \in (\pi,\theta_3 +
    \pi)$;

  \item $\theta_3 \in (\theta_2 + \pi,2\pi)$ and $\theta_4 \in (\theta_3 -
    \pi,\pi)$.
\end{itemize}
We can then compute the probability that the four bugs enter a
\texttt{cycle} state and $\theta_2 \in (0,\pi)$, by integrating
\begin{align}
  P(X_4 = \texttt{cycle} \:\cap\: \theta_2 \in (0,\pi)) =
   \frac{1}{(2\pi)^3} \int_{0}^{\pi} &\left(
    \int_{0}^{\theta_2}\int_{\pi}^{\theta_3 + \pi} 
      \, d\theta_4 \, d\theta_3 + 
    \int_{\theta_2}^{\pi}\int_{\pi}^{\theta_3 + \pi} 
      \, d\theta_4 \, d\theta_3 \,+  \right. \nonumber \\
    &\left.
    \int_{\pi}^{\theta_2 + \pi}\int_{\pi}^{\theta_3 + \pi} 
      \, d\theta_4 \, d\theta_3 + 
    \int_{\theta_2 + \pi}^{2\pi}\int_{\theta_3 - \pi}^{\pi} 
      \, d\theta_4 \, d\theta_3
  \right) \, d\theta_2 = \frac{1}{6}.
\end{align}
Doubling this result since we assumed that $\theta_2 \in (0,\pi)$, the
probability that the four bug system reaches a \texttt{cycle} steady
state is
\begin{align}
  \label{eqn:4BugsProbs}
  p_4 = \frac{1}{3}.
\end{align}
In summary, for $N \leq 4$, the probability that $N$ bugs reach a
\texttt{cycle} steady state is
\begin{align}
  \label{eqn:exactProbs}
  p_2 = 0, \quad p_3 = \frac{1}{4}, \quad p_4 = \frac{1}{3}.
\end{align}

\section{Monte Carlo Simulations for Large $N$}
\label{sec:large}
The analytic approach in section~\ref{sec:small} quickly becomes
impractical as the number of bugs, $N$, increases. Therefore, we resort
to Monte Carlo methods to estimate the probability that $N$ randomly
initialized bugs will enter the \texttt{coalesce} or \texttt{cycle}
state. In particular, we run $M$ experiments that initialize $N$ bugs
with the first bug located at $\theta_1 = 0$, while the other $N-1$ bugs
are sampled randomly from $(0,2\pi)$ with uniform distribution. Once the
bugs are initialized, their dynamics are simulated by applying Forward
Euler to the governing equation~\eqref{eqn:dynamics}. We initially used
a time step size of $\Delta t = 0.1$, but found that for larger values
of $N$, two bugs would often end up changing locations back and forth,
or a bug would pass more than one bug in a single time step, both of
which are unphysical in the continuous limit. We found that by taking
$\Delta t < \pi/N$, these unphysical dynamics are avoided, and the bugs
either enter the \texttt{cycle} or \texttt{coalesce} state. Since we
consider at most $N = 100$ bugs, we use the time step size $\Delta t =
0.01$ for all simulations.

\begin{figure}[htp]
  \centering
  \includegraphics[width=0.35\textwidth, trim=20cm 7cm 17cm 3cm,clip=true]{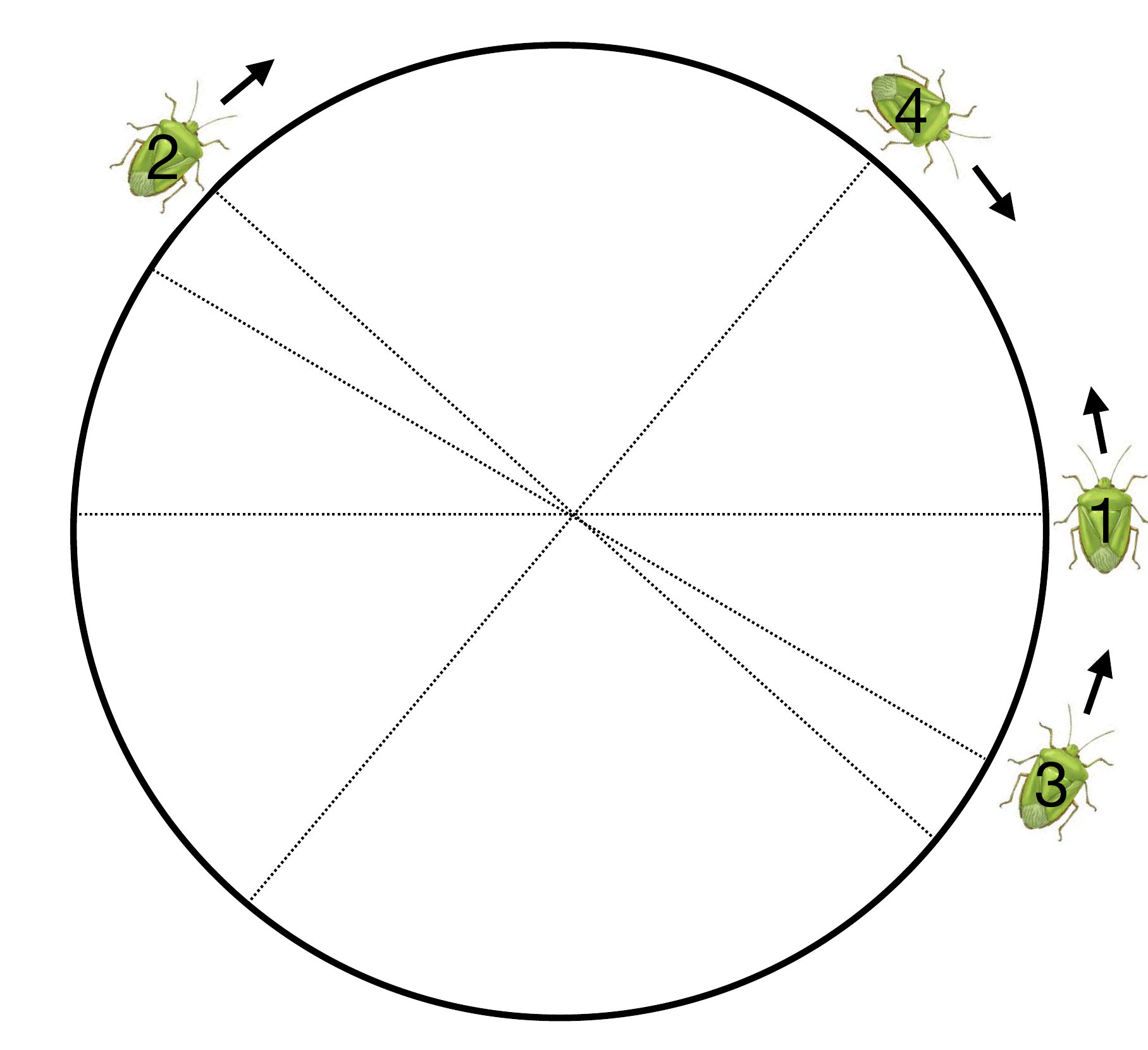}
  \caption{\label{fig:bugsCollapse} \em In this configuration, it is
  clear that the bugs will \texttt{coalesce} since all bugs lie on the
  same side of the diameter passing through bug 2 (and bug 3).}
\end{figure}

Since we only need to determine whether the bugs \texttt{coalesce} or
\texttt{cycle}, we run each simulation until we can determine which
state the system will inevitably reach. As described in the
introduction, this is accomplished by periodically checking if
\begin{itemize}
  \item all bugs are traveling in the same direction; or,
  \item all bugs are on the same side of any diameter of the circle.
\end{itemize}
The first case is easy to detect since it only requires that the
right-hand side of~\eqref{eqn:dynamics} is identically $+1$ or
identically $-1$. It is slightly harder to detect if the system has
entered the second case. In our implementation, we consider the
diameters going through each of the $N$ bugs, and if all other bugs lie
on the same side of any one of these diameters, then the bugs will
inevitably \texttt{coalesce}. Figure~\ref{fig:bugsCollapse} shows a
configuration when a simulation could be stopped since it is inevitable
that the bugs will \texttt{coalesce}.

\begin{figure}[htp]
  \centering
  \subfigimg[width=0.45\textwidth]{(a)}{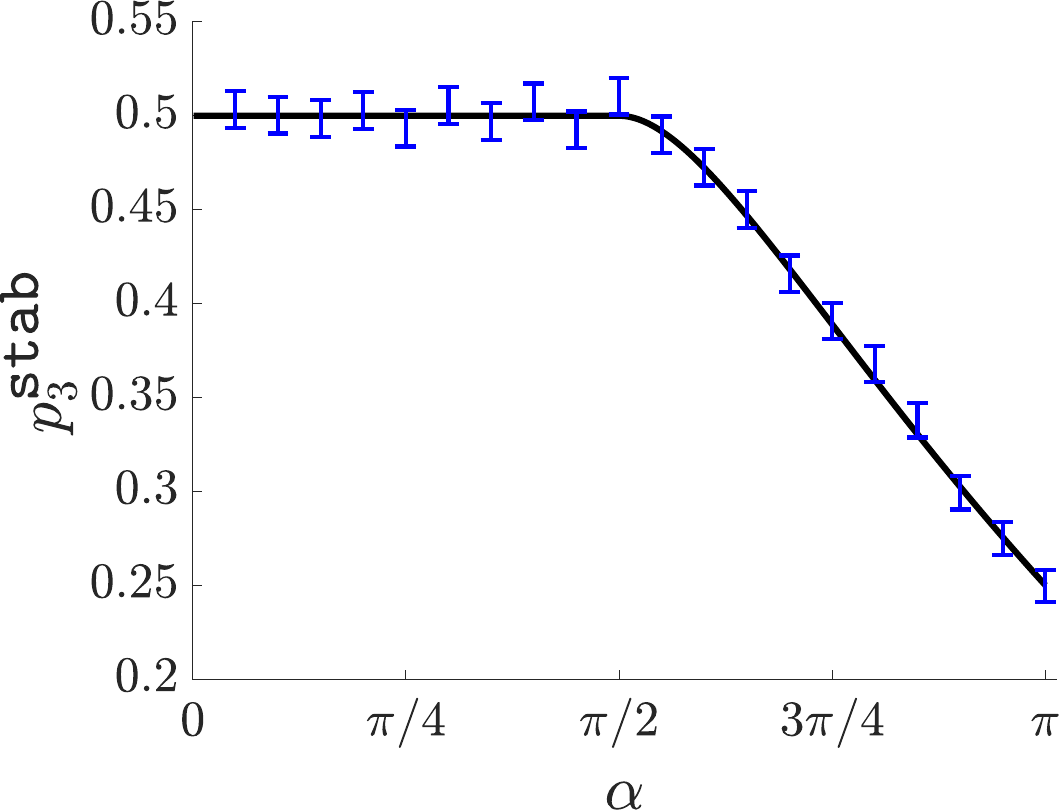}
  \quad
  \subfigimg[width=0.45\textwidth]{(b)}{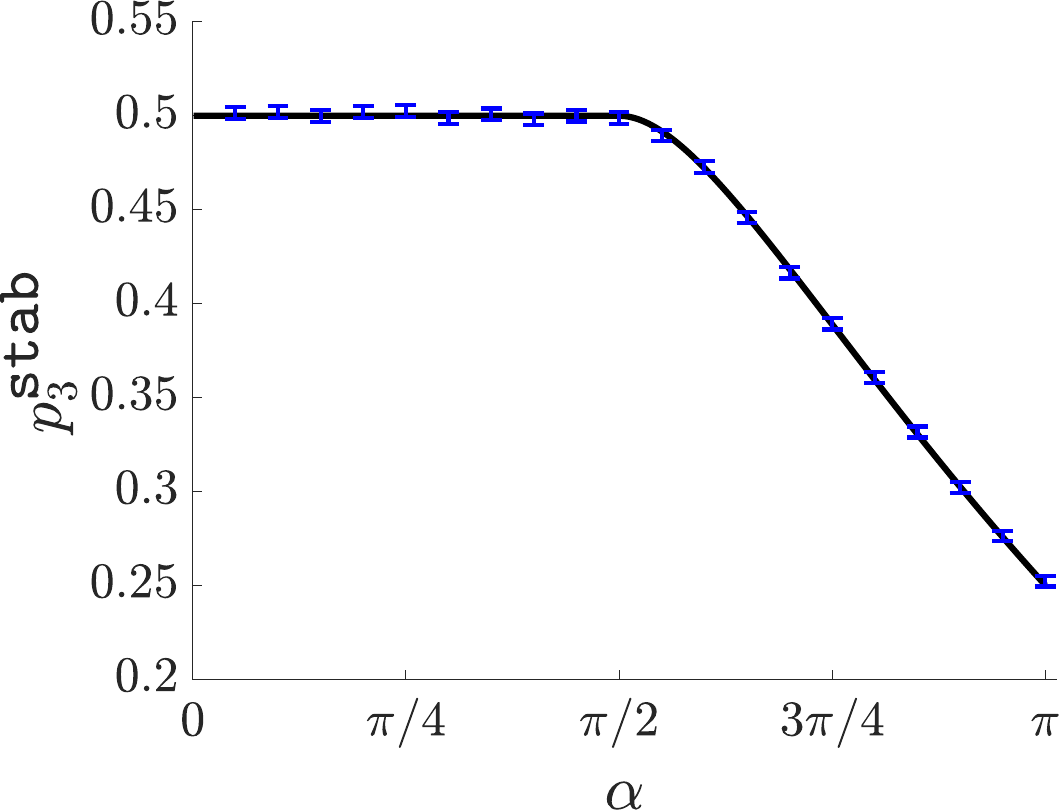}
  \caption{\label{fig:perturbation} \em Validation of the Monte Carlo
  method. The black curve is the exact probability that a perturbation
  of size no more than $\alpha$ will enter a \texttt{cycle} steady state
  (equation~\eqref{eqn:stabilityProb}). The error bars are the 95\%
  error bars of the Monte Carlo estimate using (a) $M = 10^4$
  simulations, and (b) $M = 10^{5}$ simulations.}
\end{figure}

We begin by validating our code against the analytical result for the
stability of three bugs initialized at $\theta_1 = \theta_3 = 0$ and
$\theta_2 = \pi$. This initial condition corresponds to the open black
circle in the center of Figure~\ref{fig:Stability3Bug}, and the
probability that a perturbation no larger than $\alpha$ will
\texttt{cycle} is given by equation~\eqref{eqn:stabilityProb}. We run
$M$ trials for $20$ values of $\alpha$ that are uniformly spaced in
$(0,\pi)$. When $\alpha = \pi$, perturbations of bugs 2 and 3 can lie
anywhere on the unit circle. Performing $M$ simulations provides an
estimate $\hat{p}_3^{\texttt{stab}}(\alpha)$ of
equation~\eqref{eqn:stabilityProb}. By the Central Limit Theorem, the
error in the estimate follows a normal distribution with a mean of 0 and
a variance of
\begin{align}
  \sigma^2 = \frac{\hat{p}_3^{\texttt{stab}}(\alpha)
    (1-\hat{p}_3^{\texttt{stab}}(\alpha))}{M}.
\end{align}
Figure~\ref{fig:perturbation} shows the analytical probability (solid
line) and the estimated probability using the Monte Carlo approach
(error bars) for two different values of $M$. The error bars are the
95\% confidence interval. We observe good agreement between the estimate
and exact value of $p_3^{\texttt{stab}}(\alpha)$, indicating that our
numerical approach is appropriate for estimating $p_N$ for larger $N$.

\begin{figure}[htp]
  \centering
  \subfigimg[width=0.45\textwidth]{(a)}{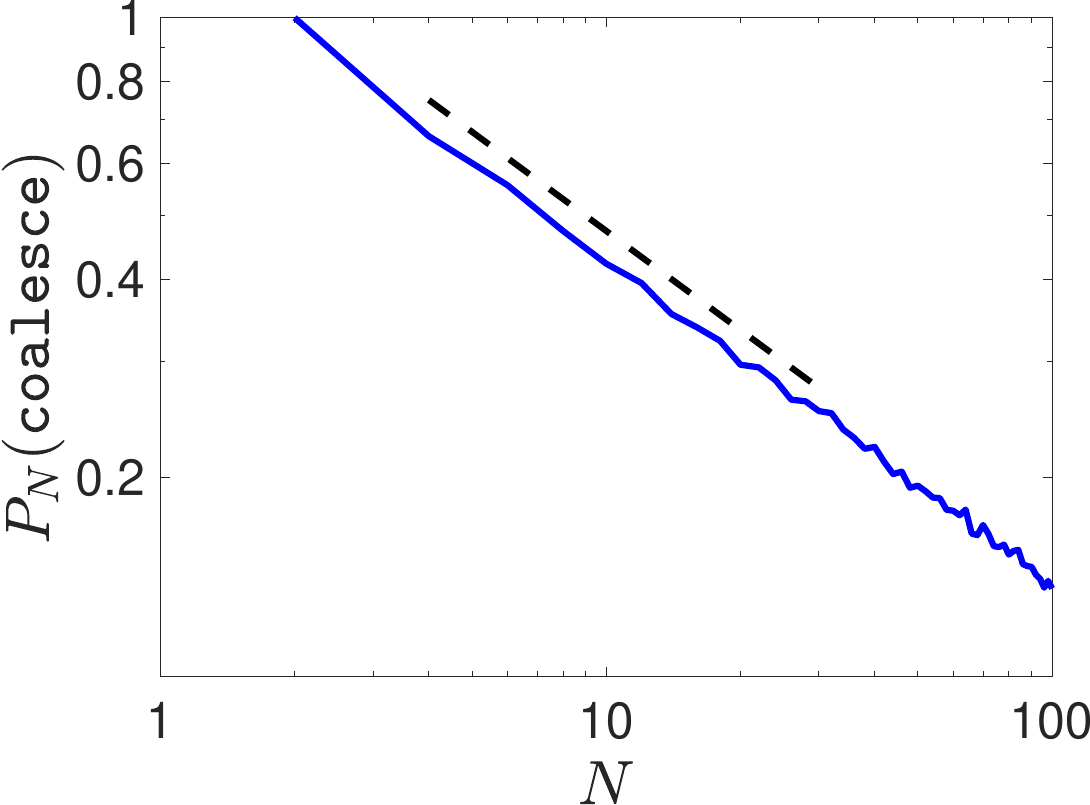}
  \quad
  \subfigimg[width=0.45\textwidth]{(b)}{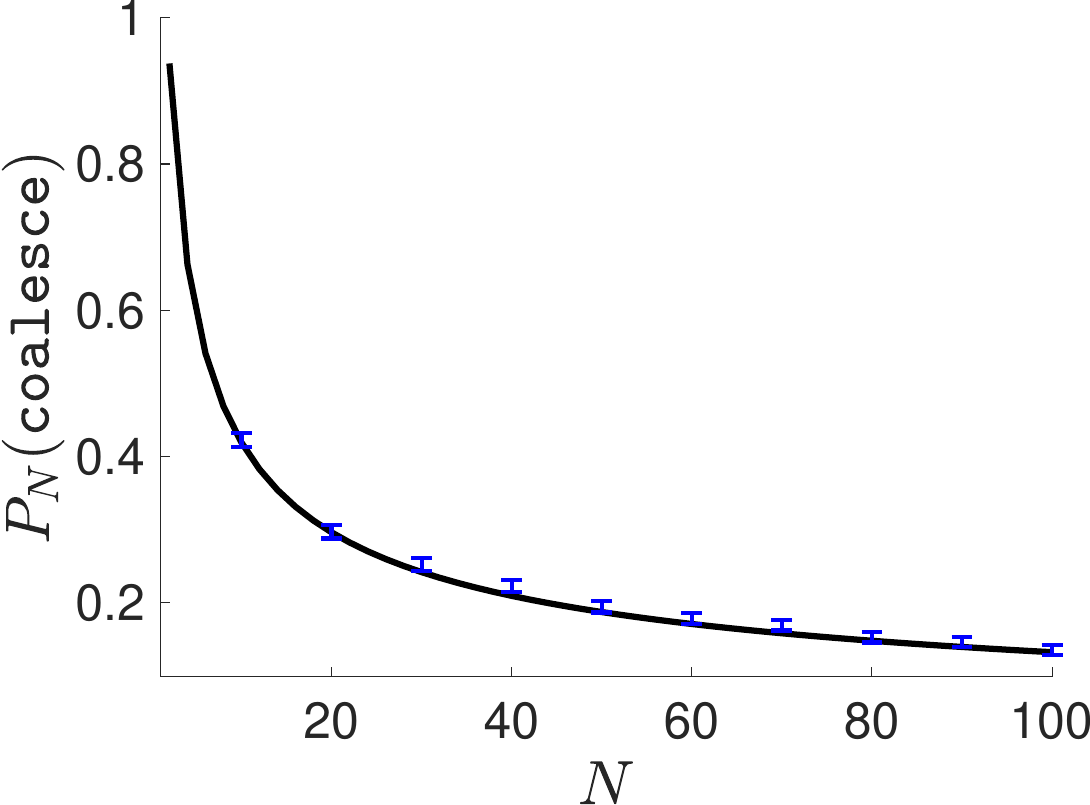}
  \caption{\label{fig:NbugCoalesce} \em (a) The estimate $\hat{p}_N$ as
  function of $N$. The black dashed line has slope $-0.5$ and indicates
  that $\hat{p}_N$ is inversely proportional to $\sqrt{N}$ (b) The least
  squares fit~\eqref{eq:coalescefunction} of the probability of the bugs
  entering a \texttt{coalesce} (solid line), and the 95\% confidence
  interval from the Monte Carlo method (intervals).}
\end{figure}

Having validated that the Monte Carlo approach converges to an analytic
result, we now estimate the probability $P(X_N = \texttt{coalesce}) = 1
- p_N$ for $N = 2,4,6,\ldots,100$. Figure~\ref{fig:NbugCoalesce} shows
the result from running $M = 10,000$ samples for each value of $N$.
First, we point out that the Monte Carlo method provides the estimates
\begin{align}
  1 - \hat{p}_2 = 1.00, \text{ and } 1 - \hat{p}^{4} = 0.661,
\end{align}
which agree with the exact values in equation~\eqref{eqn:exactProbs}.
Next, the left plot of Figure~\ref{fig:NbugCoalesce} indicates that
there is a linear relationship between $\log p_N$ and $\log N$. This
indicates a power-law relationship $P(X_N = \texttt{coalesce}) = aN^p$.
Using least squares regression, the best fit is
\begin{align}
  P(X_N = \texttt{coalesce}) \approx 1.33N^{-0.49},
  \label{eq:coalescefunction}
\end{align}
indicating that the probability of coalescing is inversely proportional
to the square root of the number of bugs. The right plot of
Figure~\ref{fig:NbugCoalesce} shows the least squares
fit~\eqref{eq:coalescefunction} (solid line) and the 95\% confidence of
the Monte Carlo method. While we see that the fit underestimates the
Monte Carlo estimate, the difference is slight, indicating that
equation~\eqref{eq:coalescefunction} provides a good estimate of the
probability that $N$ randomly initialized bugs will eventually
\texttt{coalesce}.

\section{Order Parameter and Phase Angle}
The bug configuration can be characterized using Kuramoto's order
parameter and phase angle~\cite{kuramoto1984chemical}. The
time-dependent order parameter, $r \in [0,1]$, and the phase angle,
$\psi \in [0,2\pi)$, satisfy
\begin{align}
  re^{i\psi} = \frac{1}{N}\sum_{j=1}^{N} e^{i\theta_j}.
\end{align}
A value of $r=1$ corresponds to perfect positional order, while $r=0$
corresponds to complete incoherence. If the bugs \texttt{coalesce} at
$e^{i\theta}$, then they are perfectly ordered with $r=1$, and the phase
angle is $\psi = \theta$. In contrast, if the bugs cycle clockwise
(counterclockwise), then $\psi$ decreases (increases) linearly modulo
$2\pi$, and the order parameter $r$ remains constant. The constant
equals $0$ if the cycling bugs are equally distributed, and it
approaches $1$ as the bugs begin to cluster. Therefore, the order
parameter and phase angle provide an alternative method to determine if
$N$ bugs \texttt{coalesce} or \texttt{cycle}.

\begin{figure}[htp]
  \centering
  \includegraphics[width=\textwidth]{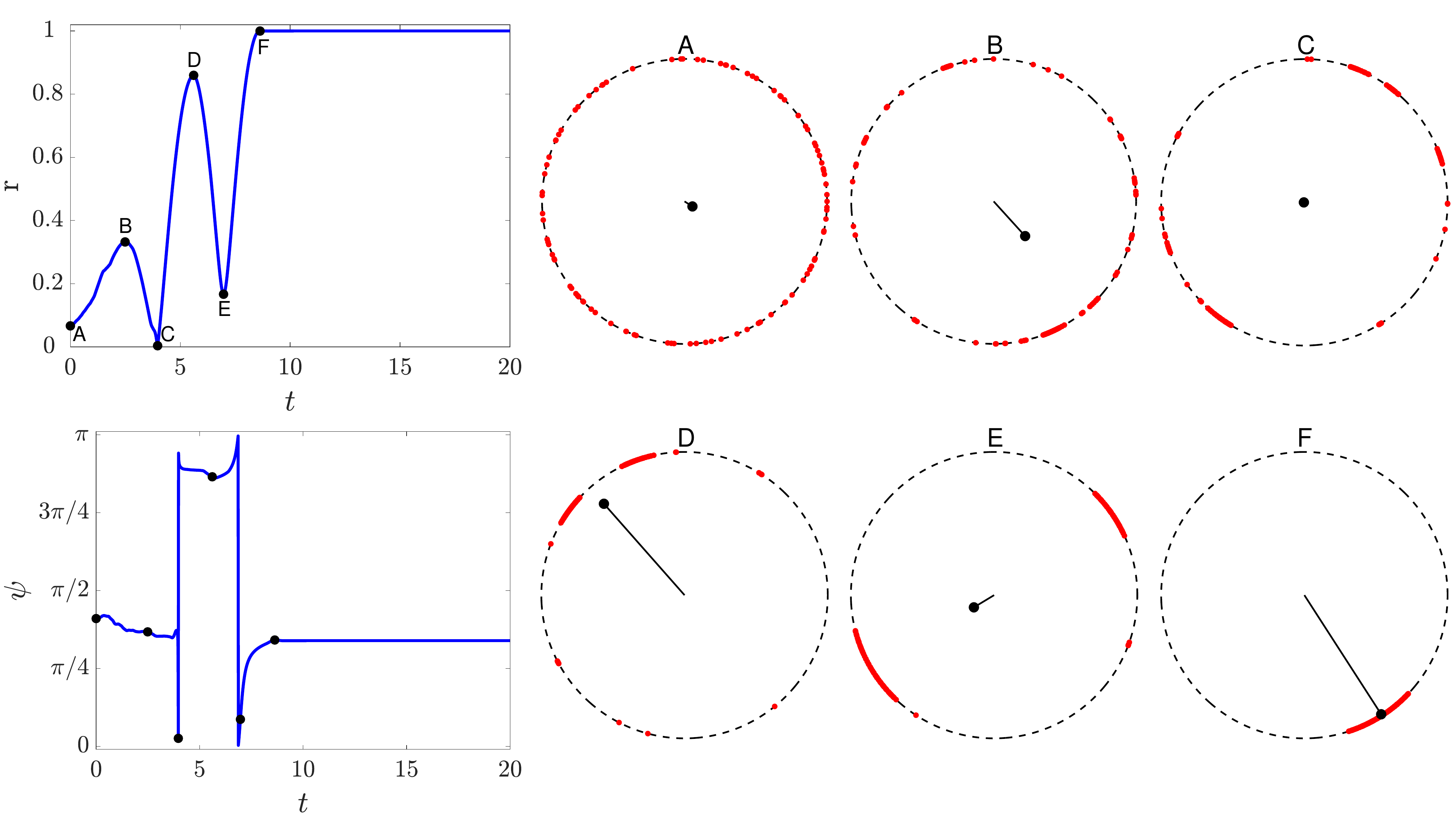}
  \caption{\label{fig:orderCoalesce} The order parameter, $r$, and phase
  angle, $\psi$, of $N=100$ bugs that eventually \texttt{coalesce}. The
  bug locations corresponding to the six black marks on the plot of the
  order parameter are shown. The vector at the center of the circle
  represents $re^{i\psi}$. Because the bugs \texttt{coalesce}, $r$
  converges to $1$ and $\psi$ converges to a constant.}
\end{figure}

\begin{figure}[htp]
  \centering
  \includegraphics[width=\textwidth]{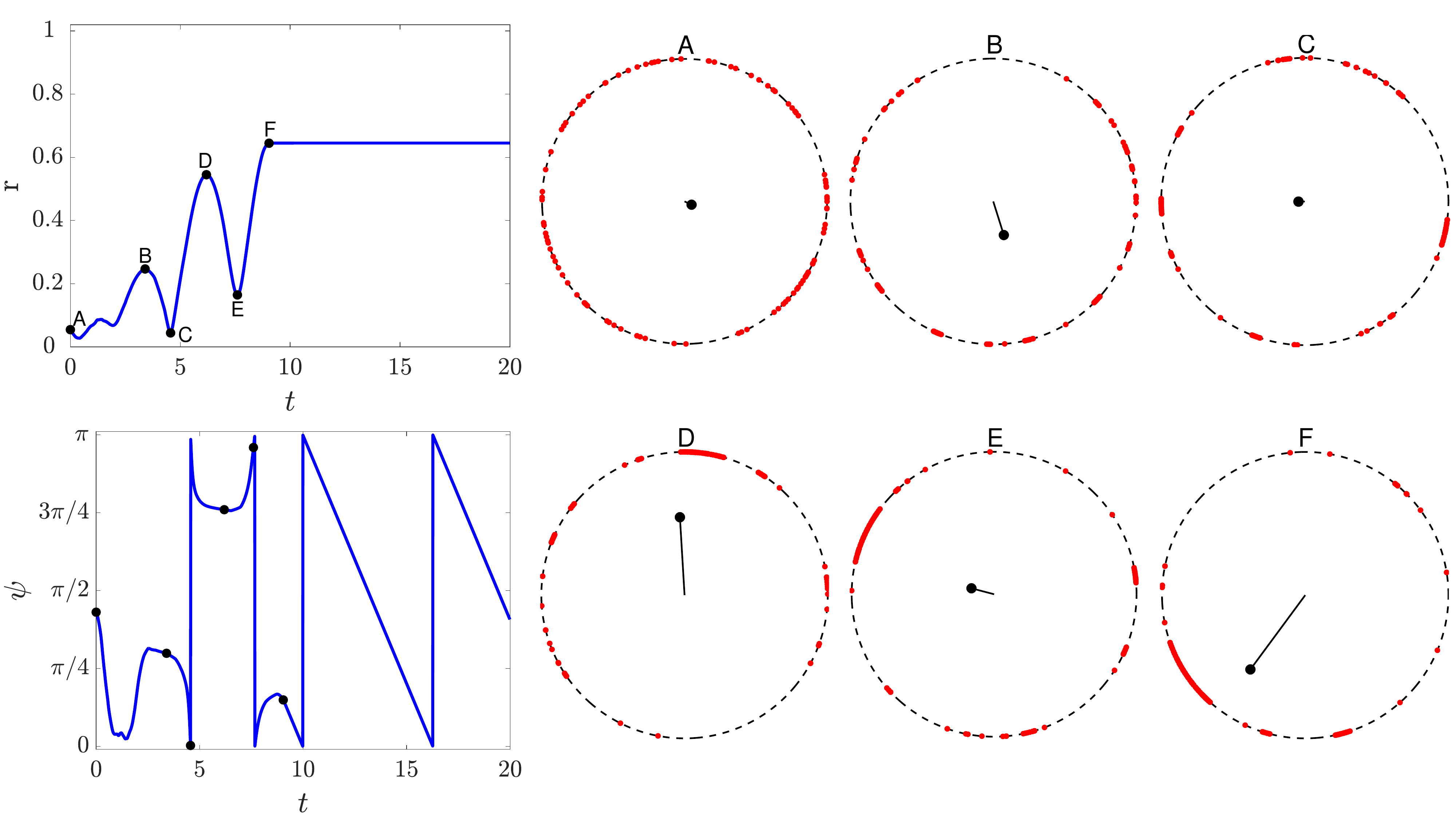}
  \caption{\label{fig:orderCycle} The order parameter, $r$, and phase
  angle, $\psi$, of $N=100$ bugs that eventually \texttt{cycle}. The bug
  locations corresponding to the six black marks on the plot of the
  order parameter are shown. The vector at the center of the circle
  represents $re^{i\psi}$. Because the bugs \texttt{cycle} clockwise,
  $r$ converges to a value less than $1$ and $\psi$ decreases linearly
  with time modulo $2\pi$.}
\end{figure}

Figure~\ref{fig:orderCoalesce} shows the order parameter and the phase
angle for $N=100$ bugs that eventually \texttt{coalesce}. The bug
locations (red dots) are shown at six critical points of the order
parameter. The values of $r$ and $\psi$ are denoted by the black line
with the black dot. We see that the order parameter can oscillate
between maximum and minimum values before ultimately reaching a value of
$r=1$, implying a \texttt{coalesce} steady state.
Figure~\ref{fig:orderCycle} shows the order parameter and phase angle of
$N=100$ bugs that eventually \texttt{cycle}. Again, the order parameter
oscillates between maximum and minimum values before ultimately reaching
a value of $r < 1$, implying a \texttt{cycle} steady state.

\section{Conclusions}

In this work, we investigated a novel generalization of the classic
cyclic pursuit problem by constraining the bugs to move exclusively
along the perimeter of a unit circle. This geometric restriction changes
the system’s long-term behavior. We identified three possible outcomes:
a stable state where all the bugs meet, defined as \texttt{coalesce}, a
stable state where all bugs infinitely chase one another, defined as
\texttt{cycle}, and an unstable state where groups of bugs are at
antipodal points of the circle, defined as \texttt{groups}. 

For small numbers of bugs ($N \leq 4$), we derived exact analytical
expressions for the probability of entering the \texttt{cycle}
steady state. These results showed that the probability of cycling
increases with $N$: $p_2 = 0$, $p_3 = 1/4$, and $p_4 = 1/3$. For larger
systems, we employed a validated Monte Carlo method to estimate the
probability that the bugs will \texttt{coalesce}. Our simulations reveal
that this probability approximately follows an inverse square-root law, 
\begin{align}
  P(X_N = \texttt{coalesce}) \approx \frac{1.33}{\sqrt{N}}.
\end{align}
The dynamics observed in this generalized setting underscore how geometric constraints can lead to qualitatively different collective behavior in pursuit problems.

\begin{figure}[htp]
  \centering
  \subfigimg[height=0.35\textwidth]{(a)}{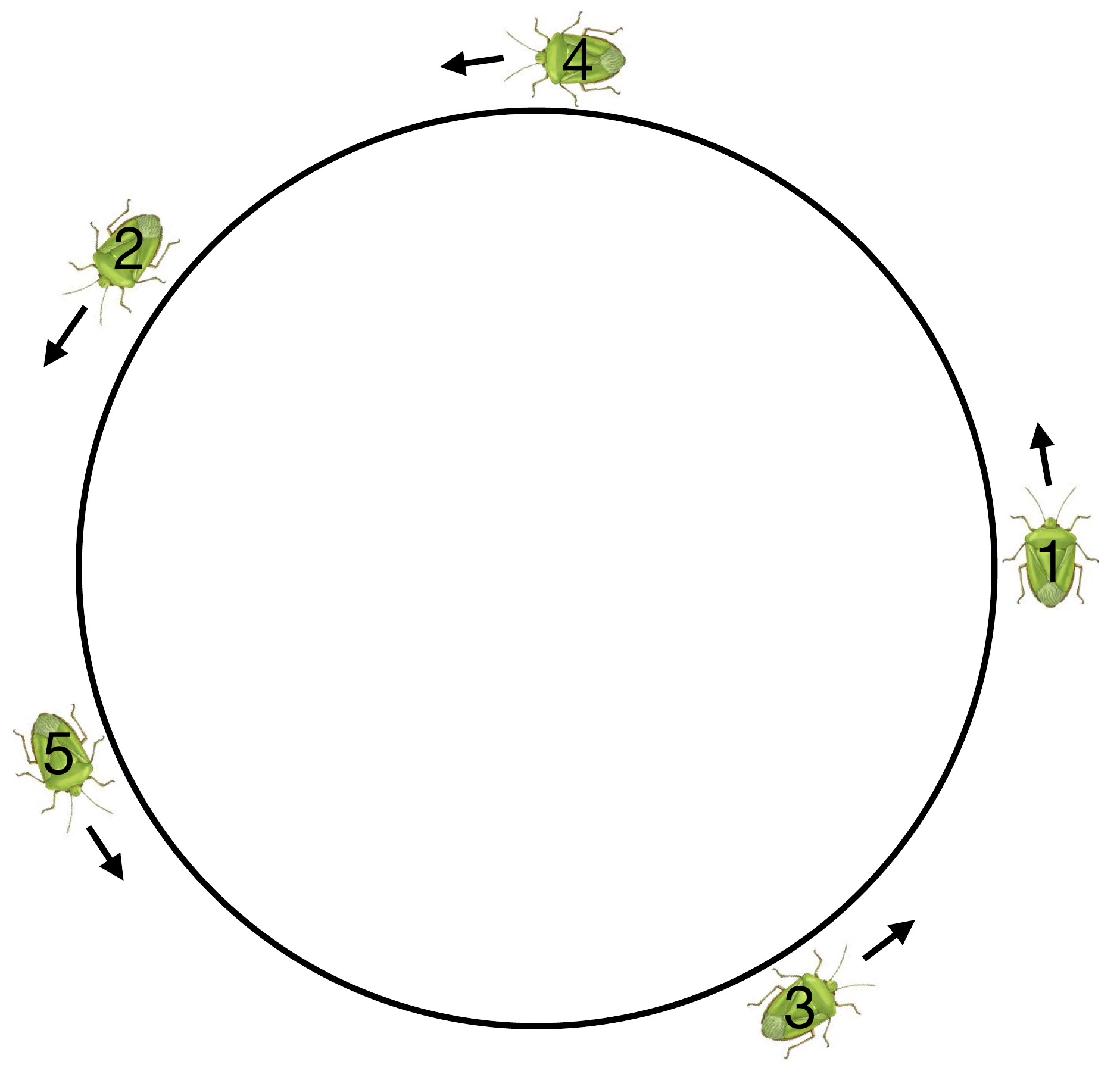}
  \quad
  \subfigimg[height=0.35\textwidth]{(b)}{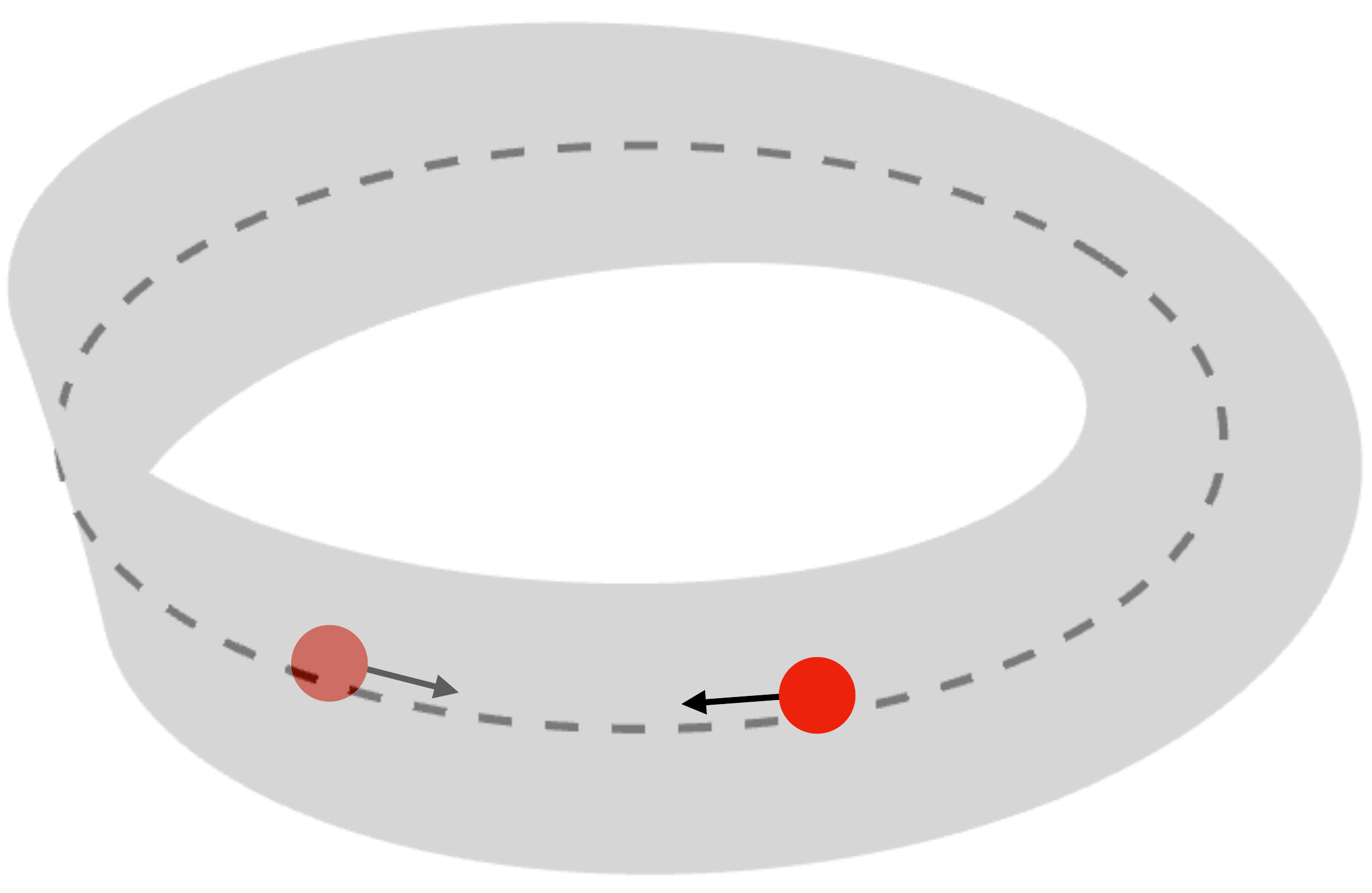}
  \caption{\label{fig:mobius} \em (a) Five bugs in a \texttt{cycle}
  steady state. Note that the indices of the bugs wrap around the circle
  twice. (b) By interpreting the M\"obius strip as a surface in
  $\mathbb{R}^3$ with thickness, it is possible for the bugs to reach a
  steady state where they are on opposite sides of from one another. In
  this figure, the bugs (red dots) are on ``opposite sides" of the
  M\"{o}bius strip.}
\end{figure}

This study establishes groundwork for several future directions. Future
work will focus on a generalization to non-uniform bug speeds. We plan
to investigate the behavior of the system if one bug is faster than the
others. We conjecture that the bugs will always \texttt{coalesce} under
this change, but this has not been explored. Additionally, there are a
variety of other statistical quantities that can be calculated
analytically for small $N$, and estimated with a Monte Carlo approach
for large $N$. For example:
\begin{itemize}
  \item What is the minimum expected wait time for $N$ randomly
  initialized bugs to enter into a stable steady state, where it is
  clear that they will either \texttt{coalesce} or \texttt{cycle} based
  on the conditions established in section~\ref{sec:large}. 

  \item For large $N$, non-consecutive bugs can cross one another many
  times before they enter into the \texttt{coalesce} or \texttt{cycle}
  steady state. What is the distribution of the number of crossings?

  \item For large $N$, the bugs can enter a \texttt{cycle} steady state
  in which each bug is moving in the same direction, but their position
  along the circle does not align with the order of their index. In this
  case, the sequence of indices wraps around the circle more than once,
  resulting in a winding number greater than one.
  Figure~\ref{fig:mobius}(a) shows $N=5$ bugs in a \texttt{cycle} steady
  state whose sequence of indices results in a winding number of two.
  For a given $N$ and an initial condition that results in the bugs
  reaching a \texttt{cycle} steady state, what is the resulting
  distribution of the winding number?
\end{itemize}

We are also considering the dynamics of $N$ randomly initialized bugs on
two-dimensional manifolds. Others have studied the dynamics of bugs
constrained to surfaces of revolution, including a
torus~\cite{sol-yeb2021}. Future directions will also consider bugs on
non-orientable surfaces. For example, by placing two bugs on the surface
of a transparent M\"obius strip embedded in $\mathbb{R}^3$, they can be
located at the same point while being on ``opposite sides" of the
surface (Figure~\ref{fig:mobius}(b)).

\paragraph{Acknowledgments:} We thank Peter Beerli and Malbor Asllani
for their valuable feedback on this paper.


\bibliographystyle{plain}

\end{document}